\documentstyle[12pt]{article}
\begin{document}

\author{S.V. Ludkowski.}

\title{$B^*$-algebras over ultranormed fields.}

\date{17 May 2018}
\maketitle

\begin{abstract}
This article is devoted to the investigation of $B^*$-algebras, dual
and annihilator ultranormed algebras. Their structure is studied in
the paper. Extensions of algebras and fields are considered and
using them core radicals and radicals are investigated. Moreover,
for this purpose also $*$-algebras and finely regular algebras are
studied. Relations with operator theory and realizations of these
algebras by operator algebras are outlined.

\footnote{key words and phrases: operator; algebra; ideal; infinite
dimension; field; ultranorm \\
Mathematics Subject Classification 2010: 12J05; 14F30; 16D60; 16E40;
46B28 \\
address: Dep. Appl. Mathematics, Moscow State Techn. Univ. MIREA,
av. Vernadksy 78, Moscow, 119454, Russia \par e-mail:
sludkowski@mail.ru }

\end{abstract}

\section{Introduction.}
Algebras and operator algebras over the real field $\bf R$ and the
complex field $\bf C$ were intensively studied. They have found
many-sided applications. For them a lot of results already was
obtained (see, for example, \cite{naribecksuffb,danschw,kadring,nai}
and references therein). Among them dual algebras and annihilator
algebras play very important role. But for such algebras over
ultranormed fields comparatively little is known because of their
specific features and additional difficulties arising from structure
of fields
\cite{ami,diarludoper,escanelb,esculbab,kersten,ludkfejms16,roo,weil}.
\par Many results in the classical case use the fact that the real
field $\bf R$ has the linear ordering compatible with its additive
and multiplicative structure and that the complex field $\bf C$ is
algebraically closed and norm complete and locally compact and is
the quadratic extension of $\bf R$, also that there are not any
other commutative fields with archimedean multiplicative norms and
complete relative to their norms besides these two fields.
\par For comparison, in the non-archimedean case the algebraic closure
of the field ${\bf Q}_p$ of $p$-adic numbers is not locally compact.
Each ultranormed field can be embedded into a larger ultranormed
field. There is not any ordering of an infinite ultranormed field
such as ${\bf Q_p}$, ${\bf C_p}$ or ${\bf F_p}(t)$ compatible with
its algebraic structure.
\par In their turn, non-archimedean analysis, functional
analysis and representations theory of groups over non-archimedean
fields develop fast in recent years
\cite{robert,roo,roosch,sch1,diar2,ludkqimb,ludkspb}. This is
motivated not only by needs of mathematics, but also their
applications in other sciences such as physics, quantum mechanics,
quantum field theory, informatics, etc. (see, for example,
\cite{anash12,aref14,aref,drag12,ish,revuzhenb,visvan2012,vla3} and
references therein).
\par This article is devoted to ultranormed $B^*$-algebras, dual algebras and annihilator
algebras over non-archimedean fields. Their structure is studied in
the paper. Extensions of algebras and fields are considered and
using them core radicals and radicals are investigated. Moreover,
for this purpose also $*$-algebras and finely regular algebras are
studied. Theorems about idempotents of algebras and their
orthogonality are proven. Division subalgebras related with
idempotents are investigated. Relations with operator theory and
realizations of these algebras by operator algebras are outlined.
Then $B^*$-algebras are defined and their properties studied.
Theorems about their embeddings into operator algebras are proved.
\par All main results of this paper are obtained for the first time.
They can be used for further studies of ultranormed algebras and
operator algebras on non-archimedean Banach spaces, their
cohomologies, spectral theory of operators, the representation
theory of groups, algebraic geometry, PDE, applications in the
sciences, etc.

\section{Ultranormed algebras and $*$-algebras.}
To avoid misunderstandings we first give our definitions and
notations.
\par {\bf 1. Notation.} Let $F$ be an infinite field supplied with a
multiplicative non-trivial ultranorm $|\cdot |_F$ relative to which
it is complete, so that $F$ is non-discrete and $\Gamma _F\subset
(0,\infty )= \{ r\in {\bf R}: 0<r<\infty \} $, where $\Gamma _F:= \{
|x|_F: x\in F\setminus \{ 0 \} \} $, whilst as usually $|x|_F=0$ if
and only of $x=0$ in $F$, also $|x+y|_F\le \max (|x|_F, |y|_F)$ and
$|xy|_F=|x|_F|y|_F$ for each $x$ and $y$ in $F$. We consider fields
with multiplicative ultranorms if something other will not be
specified.
\par If $F$ is such a field, we
denote by $E_n(F)$ the class containing $F$ and all ultranormed
field extensions $G$ of $F$ so that these $G$ are norm complete and
$|\cdot |_G|_F=|\cdot |_F$. By $E_n$ we denote the class of all
infinite non-trivially ultranormed fields $F$ which are norm
complete.
\par Henceforward, the terminology is adopted that a commutative field is
called shortly a field, while a noncommutative field is called a
skew field or a division algebra.
\par {\bf 2. Definitions.} By $c_0(\alpha ,F)$ is denoted a Banach space
consisting of all vectors $x=(x_j: ~ \forall j\in \alpha ~ x_j\in
F)$ satisfying the condition \par $card \{ j\in \alpha : ~ |x_j|>
\epsilon \} < \aleph _0$ for each $\epsilon
>0$ \\ and furnished with the norm \par $(1)$ $|x| = \sup_{j\in \alpha } |x_j|$, \\
where $\alpha $ is a set. For locally convex spaces $X$ and $Y$ over
$F$ the family of all linear continuous operators $A: X\to Y$ we
denote by $L(X,Y)$. For normed spaces $X$ and $Y$ the linear space
$L(X,Y)$ is supplied with the operator norm \par $(2)$ $|A| :=
\sup_{x\in X\setminus \{ 0 \} } |Ax|/|x|$.
\par For locally convex spaces $X$ and $Y$ over $F$
the space $L(X,Y)$ is furnished with a topology induced by a family
of semi-norms
\par $(3)$ $|A|_{p,q} :=\sup_{x\in X, p(x)>0} q(Ax)/p(x)$
\\ for all continuous semi-norms $p$ on $X$ and $q$ on $Y$.
\par Speaking about Banach spaces and Banach algebras we undermine
that a field over which it is defined is ultranorm complete.
\par If $X=c_0(\alpha ,F)$, then to each $A\in L(X,X)$ an infinite
matrix $(A_{i,j}: ~ i\in \alpha , ~ j \in \alpha )$ corresponds in
the standard basis $\{ e_j: ~ j\in \alpha \} $ of $X$, where
\par $(4)$ $x=\sum_jx_je_j$ \\
for each $x\in X=c_0(\alpha ,F)$.
\par For a subalgebra $V$ of $L(X,X)$
an operation $B\mapsto B^t$ from $V $ into $L(X,X)$ will be called a
transposition operation if it is induced by that of its infinite
matrix such that $(aA+bB)^t=aA^t+bB^t$ and $(AB)^t=B^tA^t$ and
$(A^t)^t=A$ for every $A$ and $B$ in $V$ and $a$ and $b$ in $F$,
that is $(A^t)_{i,j}=A_{j,i}$ for each $i$ and $j$ in $\alpha $.
Then $V^t := \{ A: ~ A=B^t, B\in V \} $.
\par An operator $A$ in $L(X,X)$ is called symmetric if
$A^t=A$.
\par By $L_0(X,X)$ is denoted the family of all continuous
linear operators $U: X\to X$ matrices $(U_{i,j}: ~ i\in \alpha ,
j\in \alpha )$ of which fulfill the conditions \par $(5)$ $\forall i
~\exists ~ \lim_j U_{j,i}=0$ and $\forall j ~ \exists ~ \lim_i
U_{j,i}=0$.
\par For an algebra $A$ over $F$, $F\in E_n$, it is supposed that an ultranorm
$|\cdot |_A$ on $A$ satisfies the conditions: \par $|a|_A\in (\Gamma
_F \cup \{ 0 \} )$ for each $a\in A$, also
\par $|a|_A= 0$ if and only if $a=0$ in $A$,
\par $|ta|_A=|t|_F|a|_A$ for each $a\in A$ and $t\in F$,
\par $|a+b|_A\le \max (|a|_A,|b|_A)$ and
\par $|ab|_A\le |a|_A |b|_A$ for each $a$ and $b$ in $A$. \par For
short it also will be written $|\cdot |$ instead of $|\cdot |_F$ or
$|\cdot |_A$.
\par {\bf 3. Theorem.} {\it Let $V$ be a subalgebra in $L(X,X)$ such that
$V^t=V$. Then $J$ is a left or right ideal in $V$ if and only if
$J^t$ is a right or left respectively ideal in $V$.}
\par {\bf Proof.} For each $A$ and $B$ in $V$ we get
$(AB^t)^t=BA^t$ and $B^t\in V$ and $A^t\in V$, since $V^t=V$.
Therefore, for a right ideal $J$ we deduce that $\forall A\in J ~
\forall B\in V ~ (AB^t\in J)\Leftrightarrow (BA^t\in J^t)$.
Moreover, $\forall B\in V ~ \exists U\in V ~ U^t=B$. The similar
proof is for a left ideal $J$.

\par {\bf 4. Theorem.} {\it Let $X=c_0(\alpha ,F)$, where $F\in E_n$.
Then the class $L_c(X,X)$ of all compact operators $T: X\to X$ is a
closed ideal in $L(X,X)$, also $L_{t,c}(X,X) := \{ A: ~ A\in
L_c(X,X) ~ \& ~ A^t\in L_c(X,X) \} $ is a closed ideal in
$L_0(X,X)$.}
\par {\bf Proof.} By the definition of a compact operator
$T\in L_c(X,X)$ if and only if for the closed unit ball $B$ (of
radius $1$ and with $0\in B$) in $X$ its image $TB$ is a compactoid
in $X$ (see Ch. 4 in \cite{roo}). Therefore, if $A\in L(X,X)$, then
$AB$ is bounded and convex in $X$, consequently, $TA\in L_c(X,X)$.
On the other hand, if $C$ is a compactoid in $X$, then $AC$ is a
compactoid in $X$, hence $AT\in L_c(X,X)$. Thus $L_c(X,X)$ is the
ideal in $L(X,X)$.
\par Suppose that $T_n$ is a fundamental sequence in $L_c(X,X)$
relative to the operator norm topology. Then its limit $T=\lim_n
T_n$ exists in $L(X,X)$, since $L(X,X)$ is complete relative to the
operator norm topology. Let $\epsilon >0$. There exists $m\in \bf N$
such that $|T-T_n|<\epsilon $ for each $n>m$. Since $T_n$ is the
compact operator, there exists a finite set $a_1,...,a_l$ in $X$
such that $(T_nB)\subseteq B(X,0,\epsilon )+\overline{Co}
(a_1,...,a_l)$, where \par $Co (a_1,...,a_l)=\{ x\in X: ~
x=t_1a_1+...+t_la_l, ~ t_1\in B(F,0,1,),$\par $...,t_l\in B(F,0,1)
\} $ and $B(X,y,r) := \{ z\in X: ~ |z-y|\le r \} $, \\ $~ 0<r$, $ ~
\bar{U}$ denotes the closure of a set $U$ in a topological space.
Therefore, if $x\in TB$, then there exists $y\in T_nB$ such that
$|x-y|<\epsilon $, consequently, $x\in B(X,0,\epsilon
)+\overline{Co} (a_1,...,a_l)$ due to the ultrametric inequality,
hence \par $TB \subseteq B(X,0,\epsilon )+\overline{Co}
(a_1,...,a_l)$. \par This means that the operator $T$ is compact.
Thus $L_c(X,X)$ is closed in $L(X,X)$.
\par The mapping $U\mapsto U^t$ is
continuous from $L(X,X)$ into $L(X,X)$, since $|U|=|U^t|=\sup_{i\in
\alpha , ~ j\in \alpha } |U_{i,j}|$ for each $U\in L(X,X)$. In view
of Theorem 4.39 in \cite{roo} for each $A\in L_{t,c}(X,X)$ and
$\epsilon >0$ operators $S$ and $R$ in $L(X,X)$ exist such that $SX$
and $RX$ are finite dimensional spaces over $F$ and $|A-S|<\epsilon
$ and $|A^t-R|<\epsilon $. Therefore, $L_{t,c}(X,X)\subset L_0(X,X)$
and $L_{t,c}(X,X)$ is the ideal in $L_0(X,X)$. On the other hand,
$L_0(X,X)$ is closed in $L(X,X)$, consequently, $L_{t,c}(X,X)$ is
closed in $L_0(X,X)$.

\par {\bf 5. Definition.} Suppose that $F$ is an infinite field
with a nontrivial non-archimedean norm such that $F$ is norm
complete, $F\in E_n$ and of the characteristic $char (F)\ne 2$ and
$B_2 = B_2(F)$ is the commutative associative algebra with one
generator $i_1$ such that $i_1^2=-1$ and with the involution
$(vi_1)^*=-vi_1$ for each $v\in F$. Let $A$ be a subalgebra in
$L(X,X)$ such that $A$ is also a two-sided $B_2$-module, where
$X=c_0(\alpha ,F)$ is the Banach space over $F$, $\alpha $ is a set.
We say that $A$ is a $*$-algebra if there is
\par $(1)$ a continuous bijective surjective $F$-linear operator
${\cal I}: A\to A$ such that \par $(2)$ ${\cal I}(ab)=({\cal
I}b)({\cal I}a)$ and \par $(3)$ ${\cal I}(ga)=({\cal I}a)g^*$ and
${\cal I}(ag)=g^*({\cal I}a)$ \par $(4)$ ${\cal I}{\cal I}a=a$  \par
$(5)$ $(\theta (y))(ax) = (\theta (({\cal I}a)y))(x)$ \\
for every $a$ and $b$ in $A$ and $g\in B_2$ and $x$ and $y$ in $X$,
where $\theta : X\hookrightarrow X'$ is the canonical embedding of
$X$ into the topological dual space $X'$ so that $\theta
(y)x=\sum_{j\in \alpha } y_jx_j$. For short we can write $a^*$
instead of ${\cal I}a$. The mapping $\cal I$ we call the involution.
An element $a\in A$ we call self-adjoint if $a=a^*$.

\par {\bf 6. Lemma.} {\it Let $A$ be a subalgebra of $L(X,X)$
with transposition and $A^t=A$, where $X=c_0(\alpha ,F)$, $~F\in
E_n$, $ ~ char (F)\ne 2$. Then the minimal $*$-algebra $K$ generated
by $A$ and $B_2$ has an embedding $\psi $ into $L(U,U)$ such that
$\psi (B_2)$ is contained in the center $Z(K)$ of $K$, where
$U=X\oplus X$.}
\par {\bf Proof.} We put $\psi (a) := {{a ~ 0 }\choose {0 ~ a}}$ and
$\psi (ai_1) := {{0 ~ ~a }\choose {-a ~ 0 }}$ and $(\psi (a))^* :=
{{a^t ~ 0 }\choose {0 ~ a^t}}$ and $(\psi (ai_1))^* := {{0 ~ -a^t
}\choose {a^t ~~ 0 }}$ for each $a\in A$, since $a^t\in A$.
Therefore, the minimal algebra containing $\psi (A)$ and $\psi
(Ai_1)$ is the $*$-subalgebra in $L(U,U)$. Then $(\psi
(i_1))^2=-I_U$ and $\psi (i_1) ={{0 ~ ~I_X }\choose {-I_X ~ 0 }}$,
where $I_X$ is the unit operator on $X$. Thus $\psi (ai_1)=\psi
(a)\psi (i_1)=\psi (i_1)\psi (a)=\psi (i_1a)$ for each $a\in A$ and
hence $\psi (B_2)\subset Z(K)$, where $Z(K)$ denotes the center of
the algebra $K$.

\par {\bf 7. Lemma.} {\it Let $A$ be a $*$-algebra over $F$
(see Definition 5), then each element $a\in A$ has the decomposition
$a=a_0+a_1i_1$ with $a_0^*=a_0$ and $a_1^*=a_1$ in $A$.}
\par {\bf Proof.} Put $a_0=(a+a^*)/2$, $a_1=(ai_1^*+i_1a^*)/2$, since $char (F)\ne 2$.
Then $a_0$ and $a_1$ are in $A$, since $A$ is the two-sided
$B_2$-module and $a^*\in A$ and $1\in B_2$ and $i_1\in B_2$ and
$i_1^*=-i_1$. The algebra $A$ is associative. Therefore, $a_j^*=a_j$
and $(i_1a_j)^*i_1=a_j^*i_1^*i_1=a_j=i_1(a_ji_1)^*$ for $j=0$ and
$j=1$.
\par Consider the particular case: \par if $a=a^*$, then $a_0=a$ and $(a_1i_1)^*=(a+i_1ai_1)^*/2=a_1i_1$.
\\ The latter together with $a_1^*=a_1$ implies that $-i_1a_1=a_1i_1$
if $a=a^*$. On the other hand, $a=2a_1i_1-i_1ai_1$ and $a^*=
-2a_1i_1 - i_1ai_1$ if $a=a^*$. Thus $4a_1i_1=0$ and hence $a_1=0$,
that is, $ai_1=i_1a$ if $a=a^*$, since $a_1=a_1i_1i_1^*$ and $char
(F)\ne 2$. This implies that $a_1i_1=i_1a_1$ for each $a\in A$,
consequently, the decomposition is valid $a=a_0+a_1i_1$ with the
self-adjoint elements $a_0^*=a_0$ and $a_1^*=a_1$ in $A$.

\section{Dual and annihilator ultranormed algebras.}
\par {\bf 1. Definition.} Let $A$ be a topological algebra over a field $F$
and let $S$ be a subset of $A$. The left annihilator is defined by
${\sf L}(A,S):=\{ x\in A: xS = 0 \} $ and the right annihilator is
${\sf R}(A,S):=\{ x\in A: Sx = 0 \} $, shortly they also will be
denoted by $A_l(S) := {\sf L}(A,S)$ and $A_r(S) := {\sf R}(A,S)$.
\par {\bf 2. Definition.} An
algebra $A$ is called an annihilator algebra if conditions $(1-3)$
are fulfilled:
\par $(1)$ $A_l(A)=A_r(A)=0$ and
\par $(2)$ $A_l(J_r)\ne 0$ and \par $(3)$ $A_r(J_l)\ne 0$ \\
for all closed right $J_r$ and left $J_l$ ideals in $A$.
\par If for all closed (proper or improper) left $J_l$ and right
$J_r$ ideals in $A$
\par $(4)$ $A_l(A_r(J_l))=J_l$ and
\par $(5)$ $A_r(A_l(J_r))=J_r$ \\ then $A$ is called a dual algebra.
\par If $A$ is a $*$-algebra (see Definitions 2.5) and for each $x\in A$
elements $a\in A$ and $a_1\in A$ exist such that an ultranorm on $A$
for these elements satisfies the following conditions
\par $(6)$ $|axx^*a^*_1|=|x|^2$ and $|a||a^*_1|\le 1$, \\ then the algebra $A$ is called
finely regular.

\par {\bf 3. Theorem.} {\it If $A$ is an ultranormed annihilator
finely regular Banach algebra, then $A$ is dual.}
\par {\bf Proof.} Consider arbitrary $x\in A$ and take elements $a\in A$
and $a_1\in A$ fulfilling conditions 2$(6)$, then
$|x|^2=|axx^*a^*_1|\le |a| |x| |x^*| |a^*_1|\le |x| |x^*|$, hence
$|x|\le |x^*|$. Substituting $x$ by $x^*$ we deduce analogously that
$|x^*|\le |x|$, consequently, $|x|=|x^*|$.
\par For a closed left ideal $J_l$ in $A$ if $x\in
J_l\cap (A_r(J_l))^*$, then $xx^*=0$, consequently, $x=0$ by Formula
2$(6)$ and hence $J_l\cap (A_r(J_l))^*=0$. Then $V_l := J_l\oplus
(A_r(J_l))^*$ is a left ideal in $A$, since $A_r(J_l)$ is the closed
right ideal in $A$ and $(A_r(J_l))^*$ is the closed left ideal in
$A$. \par For an arbitrary $x\in V_l$ there exist elements $y\in
J_l$ and $z\in (A_r(J_l))^*$ such that $x=y+z$. Therefore,
$xz^*=zz^*$ and $xy^*=yy^*$. Using conditions 2$(6)$ we choose
elements $a\in A$, $a_1\in A$, $b\in A$ and $b_1\in A$ with
$|a||a^*_1|\le 1$ and $|b||b^*_1|\le 1$ such that
$|azz^*a^*_1|=|z|^2$ and $|byy^*b^*_1|=|y|^2$ and hence $|x||z^*|
\ge |a||x||z^*||a^*_1|\ge |axz^*a^*_1| = |azz^*a^*_1|=|z|^2$ and
$|x| |y^*|\ge |b||x||y^*||b^*_1| \ge
|bxy^*b^*_1|=|byy^*b^*_1|=|y|^2$. Therefore, $|x|\ge |z|$ and
$|x|\ge |y|$. Thus $V_l$ is the closed left ideal in $A$.
\par From Condition 2$(3)$ it follows that a nonzero element $a\in
A$ exists such that $V_la=(0)$, consequently, $J_la=(0)$ and
$(A_r(J_l))^*a=(0)$. Then from the inclusion $a\in A_r(J_l)$ and
hence $a^*\in (A_r(J_l))^*$ it follows that $a^*a=0$. The latter
contradicts the supposition that the algebra $A$ is completely
regular. Thus $V_l=A$ and analogously for each closed right ideal
$J_r$ in $A$ the equality $A=V_r$ is valid, where $V_r= J_r\oplus
(A_l(J_r))^*$. \par Particularly, for $J_r=A_r(J_l)$ it implies that
$A=A_r(J_l)\oplus (A_l(A_r(J_l)))^*$. The involution of both sides
of the latter equality gives $A=(A_r(J_l))^*\oplus A_l(A_r(J_l))$,
since $J_l\subseteq A_l(A_r(J_l))$. Thus $J_l=A_l(A_r(J_l))$ for
each closed left ideal $J_l$ in $A$ and the involution leads to the
equality $J_r=A_r(A_l(J_r))$ for each closed right ideal $J_r$ in
$A$. Thus conditions 2$(4,5)$ are fulfilled.

\par {\bf 4. Definition.} If idempotents $w_1$ and $w_2$ of an algebra $A$ satisfy the
conditions $w_1w_2=0$ and $w_2w_1=0$, then it is said that they are
orthogonal. A family $\{ w_j: j \} $ of idempotents is said to be
orthogonal, if and only if every two distinct of them are
orthogonal. An idempotent $p$ is called irreducible, if it can not
be written as the sum of two mutually orthogonal idempotents.

\par {\bf 5. Definition.} For two Banach algebras $A$ and $B$ over an ultranormed field $F$, $F\in E_n$, we
consider  the completion $A\hat{\otimes }_{F} B$ relative to the
projective tensor product topology (see \cite{puttiel67,roo}) of the
tensor product $A\otimes _{F} B$ over the field $F$ such that
$A\hat{\otimes }_{F} B$ is also a Banach algebra into which $A$ and
$B$ have natural $F$-linear embeddings $\pi _1$ and $\pi _2$.
\par For a Banach algebra $B$ over an ultranormed field $F$, $F\in
E_n$, and an element $x\in B$ we say that $x$ has a left core
quasi-inverse $y$ if for each $H\in E_n(F)$ an element $y\in B_H$
exists satisfying the equality $x+y+yx=0$, where $B_H=B\hat{\otimes
}_FH$. Similarly is defined a right core quasi-inverse.
Particularly, if only $H=F$ is considered they are shortly called a
left quasi-inverse and a right quasi-inverse correspondingly.
\par For a unital Banach algebra $A$ over $F$, where $F\in E_n$,
if an element $x\in A$ has the property: for each field extension
$G\in E_n(F)$ the left inverse $(1+yx)^{-1}_l$ exists in $A_G$ for
each $y\in A_G$, then we call $x$ a generalized core nil-degree
element. The family of all generalized core nil-degree elements of
$A$ we call a core radical and denote it by $R_c(A)$.

\par {\bf 6. Proposition.} {\it Let $A$ be a unital Banach algebra
over $F$, where $F\in E_n$. Then \par $R_c(A) = \bigcap \{ A\cap
J_l: ~ G\in E_n(F) ~ \& ~ J_l \mbox{ is a proper maximal left ideal
in } A_G \} $.}
\par {\bf Proof.} Consider an element $x\in A$ such that for each
$G\in E_n(F)$ (see Subsection 2.1) and each maximal left ideal $J_l$
in $A_G$ the inclusion $x\in J_l$ is valid. If an element $y\in A_G$
is such that $(1+yx)^{-1}_l$ does not exist, then an element
$z=1+yx$ belongs to some left ideal $J$ in $A_G$. Since $A_G$ is the
unital algebra, then $z$ belongs to some proper maximal left ideal
$M$ such that $J\subset M$. But $yx$ also belongs to $M$, since $x$
belongs to each maximal left ideal, consequently, $1=z-yx\in M$. The
latter is impossible, since $M$ is the proper left ideal in $A_G$.
This means that the left inverse $(1+yx)^{-1}_l$ exists for every
$G\in E_n(F)$ and $y\in A_G$. Thus $x$ belongs to the core radical.
\par Vice versa. Let now $x\in R_c(A)$. Suppose the contrary that
a field extension $G\in E_n(F)$ and a proper maximal left ideal
$J_l$ in $A_G$ exist such that $x\notin J_l$. Consider the set $V$
of all elements $z=b-yx$ with $b\in J_l$ and $y\in A_G$. Evidently
$V$ is the left ideal in $A_G$ containing $J_l$, but $J_l$ is
maximal, consequently, $V=A_G$. This implies that $1=b-yx$ for some
$b\in J_l$ and $y\in A_G$. Therefore the element $b=1+yx$ has not a
left inverse. But this contradicts the supposition made above.

\par {\bf 7. Proposition.} {\it Suppose that $A$ is a unital Banach algebra
over $F$, where $F\in E_n$. Then \par $(x\in R_c(A)) \Leftrightarrow
(\forall G\in E_n(F) ~ \forall y\in A_G ~ \exists (1+yx)^{-1} \in
A_G )$.}
\par {\bf Proof.} If $\forall G\in E_n(F) ~ \forall y\in A_G ~ \exists (1+yx)^{-1} \in
A_G$, then $\forall G\in E_n(F) ~ \forall y\in A_G ~ \exists
(1+yx)^{-1}_l \in A_G$, consequently, $x\in R_c(A)$, where as
usually $(1+yx)^{-1}$ notates the inverse of $1+yx$.
\par Vice versa. Let $x\in R_c(A)$. Then by the definition
of the core radical $\forall G\in E_n(F) ~ \forall y\in A_G ~
\exists (1+yx)^{-1}_l \in A_G$. For $G\in E_n(F)$ denote by $1+b$ a
left inverse of $1+yx$ in $A_G$, that is $(1+b)(1+yx)=1$. This
implies that $1+yx$ is the right inverse of $1+b$ in $A_G$ and
$b=-byx-yx$. From $x\in R_c(A)$ it follows that $b\in R_c(A_G)$,
since $x\in J_l$ and hence $y\in J_l$ for each proper maximal left
ideal $J_l$ in $A_H$ and each $H\in E_n(G)$. This means that for
every $H\in E_n(G)$ and $z\in A_H$ a left inverse $(1+zb)^{-1}_l$
exists in $A_H$, particularly, for $z=1$ also. On the other hand,
the right inverse is $(1+zb)^{-1}_r=1+yx$ as it was already proved
above. Therefore the inverse (i.e. left and right simultaneously)
$(1+b)^{-1}=1+yx$ exists. Thus $1+b$ is the inverse of $1+yx$ in
$A_G$.

\par {\bf 8. Proposition.} {\it Let $A$ be a unital Banach algebra
over $F$, where $F\in E_n$. Then \par $R_c(A) = \bigcap \{ A\cap
J_r: ~ G\in E_n(F) ~ \& ~ J_r \mbox{ is a proper maximal right ideal
in } A_G \} $. Moreover, $R_c(A)$ is the two-sided ideal in $A$.}

\par {\bf Proof.} Consider the class $Q_e(A)$ of all elements $x\in A$ such that
for each field extension $G\in E_n(F)$ the right inverse
$(1+xy)^{-1}_r$ exists in $A_G$ for each $y\in A_G$. Analogously to
the proof of Proposition 6 we infer that
\par $Q_e(A) = \bigcap \{
J_r: ~ G\in E_n(F) ~ \& ~ J_r \mbox{ is a proper maximal right ideal
in } A_G \} $. \\ Similarly to the proof of Proposition 7 we deduce
that
\par $(x\in Q_e(A)) \Leftrightarrow (\forall G\in E_n(F) ~
\forall y\in A_G ~ \exists (1+xy)^{-1} \in A_G )$.
\par Suppose that $G\in E_n(F)$, $x\in A$, $y\in A_G$ and the
inverse  element exists $(1+yx)^{-1}=1+b$ in $A_G$. Then
$(1+xy)(1-xy-xby)-1=-x((1+yx)(1+b)-1)y=0$ and $(1-xy-xby)(1+xy)-1=
-x((1+b)(1+yx)-1)y=0$, consequently, $1-xy-xby =(1+xy)^{-1}$.
Analogously if the inverse element $(1+yx)^{-1}$ exists, then
$(1+xy)^{-1}$ also exists. This implies that $Q_e(A)=R_c(A)$ and
hence the core radical is the two-sided ideal in $A$.

\par {\bf 9. Proposition.} {\it Suppose that $A$ is a unital Banach algebra
over $F$, where $F\in E_n$. Then an extension field $H=H_F\in
E_n(F)$ exists such that $R_c(A)=A\cap R(A_H)$, where $R(A_H)$
denotes the radical of the algebra $A_H$ over $H$. Moreover, $H$ can
be chosen algebraically closed and spherically complete.}
\par {\bf Proof.} Consider an arbitrary element $x\in A\setminus
R_c( A)$. This means that a field extension $G=G_x\in E_n(F)$ and an
element $y\in A_G$ exist such that the element $(1+yx)$ has not the
left inverse in $A_G$. For the family $ {\cal G} := \{ G_x: ~ x\in
A\setminus R_c( A), G_x\in E_n(F) \} $ a field $H=H_F\in E_n(F)$
exists such that $G_x\subseteq H$ for each $x\in A\setminus R_c( A)$
due to Proposition V.3.2.2 \cite{bourbcorcb} and since the
multiplicative ultranorm $|\cdot |_F$ can be extended to a
multiplicative ultranorm $|\cdot |_H$ on $H$ (see Proposition 5 in
Section VI.3.3 \cite{bourbcomalgb}, Krull's existence theorem 14.1
and Theorem 14.2 in \cite{sch1} or 3.19 in \cite{roo}, Lemma 1 and
Proposition 1 in \cite{diar})). \par If $H_F$ is not either
algebraically closed or spherically complete, one can take the
spherical completion of its algebraic closure $\bar{H}_F$ (see
Corollary 3.25, Theorem 4.48 and Corollary 4.51 in \cite{roo}). Then
also $\bar{H}_F\in E_n(F)$. Denote shortly $\bar{H}_F$ by $H$.
\par Therefore, if $G\in \cal G$, then from $y\in A_G$ it follows
that $y\in A_H$. For each $x\in A\setminus R_c( A)$ an element $y\in
A_G$ exists such that $A_G(1+yx)$ is a left proper ideal in $A_G$,
consequently, $A_G(1+yx)\hat{\otimes }_GH= A_H(1+yx)$ is a left
proper ideal in $A_H$, since $H\subset Z(A_H)$. Therefore, $(1+yx)$
has not a left inverse in $A_H$.
\par Thus for each $x\in A\setminus R_c(A)$ and $G=G_x\in \cal
G$ and element $y\in A_H$ exists such that $(1+yx)$ has not a left
inverse in $A_H$. Therefore, $A\cap R(A_H) \subset R_c(A)$. On the
other hand, if $x\in R_c(A)$, then $x\in R(A_H)$ according to the
definition of $R_c(A)$ in \S 5. Thus $R_c(A)=A\cap R(A_H)$ for the
fields $H$ constructed above.

\par {\bf 10. Theorem.} {\it Let $A$ be a unital Banach algebra
over $F$, where $F\in E_n$. Then an extension field $K=K_F\in
E_n(F)$ exists such that \par $(1)$ $R_c(A_K)=R(A_K)$. \par
Moreover, $K$ can be chosen algebraically closed and spherically
complete.}
\par {\bf Proof.} Put $K_1=H$, where $H=H_F$ is given by Proposition 9.
Then by induction take $K_{n+1}=H_{K_n}$ for each natural number
$n=1, 2, 3,...$. There are isometric embeddings $K_n\hookrightarrow
K_{n+1}$ for each $n$. Let $K$ be the norm completion of $K_{\infty
} := \bigcup_{n=1}^{\infty } K_n$, hence $K\in E_n(F)$. In addition
each field $K_l$ can be chosen algebraically closed and spherically
complete due to Proposition 9. Moreover, it is possible to take as
$K$ the spherical completion of the algebraic closure of $K_{\infty
}$ (see Corollary 3.25, Theorem 4.48 and Corollary 4.51 in
\cite{roo}).
\par In view of Proposition 9 $R_c(A_{K_l})=A_{K_l}\cap
R(A_{K_{l+1}})$ for each natural number $l$. Let $x\in R_c(A_K)$,
that is for each $G\in E_n(K)$ and $y\in A_K$ a left inverse
$(1+yx)_l^{-1}$ exists in $A_K$. The algebra $A\otimes _F{K_{\infty
}}$ over the field $K_{\infty }$ is everywhere dense in
$A_K=A\hat{\otimes }_FK$. Therefore, there exist sequences $x_n$ and
$y_n$ in $A_K$ such that $x_n\in A_{K_n}$ and $y_n\in A_{K_n}$ for
each $n$ and $\lim_n x_n=x$ and $\lim_n y_n=y$. Since $(1+z)$ is
invertible in $A_K$ for each $z\in A_K$ with $|z|<1$, then a natural
number $m$ exits such that a left inverse $(1+y_nx_n)_l^{-1}$ exists
for each $n>m$. \par From $G\in E_n(K)$ and $K_l\in E_n(F)$,
$K_l\subseteq K$ it follows that $G\in E_n(K_l)$ for each $l=1, 2,
3,...$.
 On the other hand,
an element $y\in A_K$ can be any marked element in particularly
belonging to $A_{K_l}$. Thus $\bigcup_l R_c(A_{K_l})$ is dense in
$R_c(A_K)$. Similarly considering $G=K$ one gets that $\bigcup_l
R(A_{K_l})$ is dense in $R(A_K)$. Mentioning that $\bigcup_l
A_{K_l}$ is dense in $A\otimes _F{K_{\infty }}$ one gets that
$\bigcup_l A_{K_l}$ is dense in $A_K$. Therefore, we infer that
$$R_c(A_K)=cl_{A_K}(\bigcup_l R_c(A_{K_l})) = cl_{A_K}(\bigcup_l
(A_{K_l}\cap R(A_{K_{l+1}})))$$  $$ = cl_{A_K}(\bigcup_l
R(A_{K_{l+1}}))= R(A_K),$$  where $cl_{A_K} B$ denotes the closure
of a subset $B$, $B\subset A_K$, in $A_K$.

\par {\bf 11. Proposition.} {\it Let $A$ be a Banach
algebra over $F$, $F\in E_n$, also let a field $K$ fulfill Condition
10$(1)$ for $A_1$, where $A_1=A$ if $1\in A$, while $A_1=A\oplus 1F$
if $1\notin A$. Then an element $x\in A_K$ is not core left
quasi-invertible if and only if $J_{l,G} := \{ z+zx: ~ z\in A_G \} $
is a proper left ideal in $A_G$ for each $G\in E_n(K)$. If so
$J_{l,G}$ is a proper regular left ideal in $A_G$ such that $x\notin
J_{l,G}$.}
\par {\bf Proof.} By virtue of Theorem 10 $R_c(A_{1,K})=R(A_{1,K})$.
Hence for each $G\in E_n(K)$ an element $x\in A_K$ is not core left
quasi-invertible in $A_G$ if and only if it does not belong to
$R(A_{1,K})$. If $u=y+yx$ and $v=z+zx$ belong to $J_{l,G}$, $b$ and
$c$ are in $G$, where $y$ and $z$ belong to $A_G$, then
$bu+cv=(by+cz)+(by+cz)x$, consequently, $cu+bv\in J_{l,G}$. That is
$A_GJ_{l,G}\subseteq J_{l,G}$. If $J_{l,G}$ is not a proper left
ideal, then $J_{l,G}=A_G$. This implies that an element $z\in A_G$
exists such that $x+zx=-x$. The latter is equivalent to the equality
$x+zx+x=0$. Thus $z$ is a left quasi-inverse of $x$.
\par Vise versa if $x$ has a left quasi-inverse in $A_G$, then $x\in
J_{l,G}$, hence $-zx\in J_{l,G}$. Therefore, $z=(z+zx)-zx\in
J_{l,G}$ for each $z\in A_G$, consequently, $A_G=J_{l,G}$. Thus if
$J_{l,G}$ is a proper left ideal, then $x\notin J_{l,G}$. Mention
that the element $w=-x$ is unital modulo the proper left ideal
$J_{l,K}$, consequently, this ideal is regular.

\par {\bf 12. Proposition.} {\it Suppose that $A$ is a Banach
algebra over $F$, $F\in E_n$, also a field $K$ satisfies Condition
10$(1)$ for $A_1$. Then the following conditions are equivalent:
\par $(1)$ an element
$x\in A_K$ possesses a left quasi-inverse in $A_G$ for each $G\in
E_n(K)$;
\par $(2)$ for every $G\in E_n(K)$ and a maximal regular
proper left ideal $M_{l,G}$ in $A_G$ an element $y\in A_G$ exists
such that $x+y+yx\in M_{l,G}$.}
\par {\bf Proof.} If an element $x\in A_K$ possesses a left quasi-inverse $y_G$ in $A_G$
for each $G\in E_n(K)$, then $x+y_G+y_Gx\in M_{l,G}$ for each
maximal regular proper left ideal $M_{l,G}$ in $A_G$ due to Theorem
10.
\par Vise versa suppose that Condition $(2)$ is fulfilled, but $x$ is not
left quasi-invertible in $A_G$ for some $G\in E_n(K)$. Then
$J_{l,G}$ is a regular proper left ideal in $A_G$ according to
Proposition 11. Therefore, a maximal regular proper left ideal
$M_{l,G}$ in $A_G$ exists containing $J_{l,G}$. Thus an element
$y\in A_G$ exists such that $x+y+yx\in M_{l,G}$. On the other hand,
the inclusion $y+yx\in J_{l,G}$ is accomplished, consequently, $y\in
M_{l,G}$ and hence $-zx\in M_{l,G}$ for each $z\in A_G$. This
implies that $z\in M_{l,G}$ for each $z\in A_G$, since
$z=-zx+(z+zx)$. But this leads to the contradiction $A_G=M_{l,G}$.
Thus $(2)\Rightarrow (1)$.

\par {\bf 13. Proposition.} {\it Suppose that $A$ is a Banach annihilator algebra
over an ultranormed field $F$, $F\in E_n$. Then a field extension
$K$, $K\in E_n(F)$, exists such that if an element $-p\in A_K$ is
not core left quasi-invertible, then a nonzero element $x\in
A_K\setminus \{ 0 \} $ exist satisfying the equation $px=x$.}
\par {\bf Proof.} We take a field $K$, $K\in E_n(F)$, given by Theorem 10 for a unital
algebra $E=A_1$, where $E=A\oplus 1F$ if $1\notin A$, while $E=A$ if
$1\in A$. Therefore, $R_c(E_K)=R(E_K)$.
\par By virtue of Proposition 11 $J_{l,K} := \{ yp-p: y\in A_K \} $
is a regular proper left ideal in $A_K$. Since $E_K$ is the unital
Banach algebra over $K$, then it is with continuous inverse. Hence
if $A$ is not unital, then $A_K$ is with the continuous
quasi-inverse. Mention that an element $v$  is a left quasi-inverse
of $q$ in $A_K$ if and only if $1+v$ is a left inverse of $1+q$ in
$E_K$.
\par Therefore, if $1\notin A$, then a bijective correspondence exists:
$Q$ is a left (maximal) ideal of $E_K$ which is not contained
entirely in $A$ if and only if $Q\cap A_K$ is a regular (maximal
respectively) left ideal of $A_K$. If $1\in A$, then each left ideal
in $A_K$ is regular.
\par Recall that a ring $B$ satisfying the identities \par $(1)$ ${\sf L}(B,B)=(0)$ and
${\sf R}(B,B)=(0)$ is called annihilator, where \par $(2)$ ${\sf
L}(B,S) =
\{ x\in B: xS=(0) \} $ and ${\sf R}(B,S) = \{ x\in B: Sx=(0) \} $ \\
denote a left annihilator and a right annihilator correspondingly of
a subset $S$ in $B$. Thus \par $(3)$ ${\sf L}(A_K,A_K)=(0)$ and
${\sf R}(A_K,A_K)=(0)$,
\\ since $A_K=A\hat{\otimes }_FK$, since by the conditions of this
proposition ${\sf L}(A,A)=(0)$ and ${\sf R}(A,A)=(0)$, also $A$ and
$A_K$ are Banach algebras. Next we take the closure $cl_{A_K}
(J_{l,K})$ of $J_{l,K}$ in $A_K$. Therefore, ${\sf R}(A_K,cl_{A_K}
(J_{l,K}))$ is not nil, ${\sf R}(A_K,cl_{A_K} (J_{l,K}))\ne (0)$.
\par Suppose that $x$ is a nonzero element in ${\sf R}(A_K,cl_{A_K}
(J_{l,K}))$, consequently, $x\in {\sf R}(A_K,J_{l,K})$.
\par If $z\in {\sf R}(A_K,J_{l,K})$, then $y(pz-z)=(yp-y)z=0$ for each
$y\in A_K$. From ${\sf L}(A_K,A_K)=(0)$ and ${\sf R}(A_K,A_K)=(0)$
it follows that $pz-z=0$. Vise versa, if $pz-z=0$ for some $z\in
A_K$, then $(yp-y)z=y(pz-z)=0$ and hence $z\in {\sf
R}(A_K,J_{l,K})$. Therefore,
\par $(4)$ ${\sf R}(A_K,J_{l,K})=\{ z\in A_K: pz=z \} $. \par Thus $px=x$.

\par {\bf 14. Theorem.} {\it Suppose that $A$ is a Banach annihilator algebra
over a field $F\in E_n$ such that $R_c(A)=R(A)$ and $M_r$ is a
proper maximal closed right ideal in $A$ satisfying the condition
${\sf L}(A,M_r)\cap R(A)=(0)$. Then ${\sf L}(A,M_r)$ contains an
idempotent $p$ and \par $(1)$ ${\sf L}(A,M_r)=Ap$ and \par $(2)$
$M_r= \{ z-pz: z\in A \} $.}
\par {\bf Proof.} A nonzero element $b$ in ${\sf L}(A,M_r)$ exists, since
${\sf L}(A,M_r)\ne (0)$, since $M_r$ is a proper right ideal in $A$.
Therefore, $M_r\subset {\sf R}(A, \{ b \}  )\ne A$ and consequently,
\par $(3)$ ${\sf R}(A, \{ b \}  ) = M_r$, \\ since the right ideal $M_r$ is maximal.
The element $b$ does not belong to $R(A)$, since ${\sf L}(A,M_r)\cap
R(A)=(0)$ by the conditions of this theorem. \par In view of Theorem
10 and Propositions 11 and 12 a scalar $t\in F$ and an element $y\in
A$ exist such that the element $-p=tb+yb$ has not a left
quasi-inverse in $A_G$ for each $G\in E_n(F)$. Thus $p\ne 0$ and
$p\in {\sf L}(A,M_r)$. By virtue of Proposition 13 a nonzero element
$x\in A\setminus (0)$ exists such that $px=x$, consequently,
$(p^2-p)x=0$.
\par Suppose that $p^2-p$ is not nil, $p^2-p\ne 0$. We have
$p^2-p\in {\sf L}(A,M_r)$. Taking $b=p^2-p$ in $(3)$ one gets ${\sf
R}(A,p^2-p)=M_r$, consequently, $(p^2-p)x\in M_r$ and inevitably
$x=px=0$. This leads to the contradiction. Thus $p^2=p$.
\par On the other hand, $p\in {\sf L}(A,M_r)$ and $p$ is not nil.
Taking $b=p$ in $(3)$ provides $M_r={\sf R}(A,\{ p \} )$ and ${\sf
R}(A,\{ p \} ) = \{ z-pz: z\in A \} $, since $p(y-py)=py-p^2y=0$,
also if $pz=0$, then $z=z-pz$.  Therefore, ${\sf L}(A,M_r)=Ap$ due
to 13$(4)$ and since $p$ is the idempotent.

\par {\bf 15. Corollary.} {\it If conditions of Theorem 14 are
fulfilled, then $M_r$ is a maximal right ideal and ${\sf L}(A,M_r)$
is a minimal left ideal, also $pA$ is a minimal right ideal and
${\sf L}(A,pA)$ is a maximal left ideal.}

\par {\bf 16. Theorem.} {\it Let $A$ be a Banach annihilator algebra
over a field $F\in E_n$ such that $R_c(A)=R(A)$, let also $J_l$ be a
minimal left (may be closed) ideal which is not contained in $R(A)$,
$J_l\setminus R(A)\ne \emptyset $. Then $J_l$ contains an idempotent
$p$ for which $J_l=Ap$ and ${\sf R}(A,Ap) = \{ x-px: x\in A \} $.}
\par {\bf Proof.} Take $x\in J_l\setminus R(A)$. From Propositions
11 and 12 it follows that $b\in F$ and $y\in A$ exist such that the
element $-p=bx+yx$ has not a left quasi-inverse, consequently, $p\ne
0$. \par In view of Proposition 13 an element $v\in A$ exists having
the property $pv=v$. Therefore, $Y_l := \{ z\in J_l: xv=0 \} $ is a
left ideal such that it is contained in $J_l$ and $J_l\ne Y_l$,
since $p\in J_l\setminus Y_l$. This ideal $Y_l$ is closed, if $J_l$
is closed. The ideal $J_l$ is minimal, hence $Y_l=(0)$. This implies
that $zv\ne 0$ if $z\in J_l\setminus \{ 0 \} $.On the other hand,
$p^2-p\in J_l$ and $(p^2-p)v=0$, hence $p^2-p=0$. Thus $p$ is the
idempotent. \par For each $z\in Ap$ the condition $z=zp$ is valid,
consequently, $Ap$ is a closed left ideal contained in $J_l$ and
hence $Ap=J_l$, since the left ideal $J_l$ is minimal. Therefore,
${\sf R}(A,J_l) = \{ x-px: x\in A \} $.

\par {\bf 17. Lemma.} {\it If $A$ is a Banach annihilator semi-simple algebra
over a field $F\in E_n$ with $R_c(A)=R(A)$ and $J$ is a left (or
right, or two-sided) ideal in $A$ such that $J^2=(0)$, then
$J=(0)$.}
\par {\bf Proof.} Suppose that $J$ is a left ideal in $A$ with
$J^2=(0)$. Therefore, $(tx+yx)^2=0$ for every $t\in F$, $x\in  J$
and $y\in A$, since $tx+yx\in J$. In this case the element $z=tx+yx$
has the left quasi-inverse $-z$. By virtue of Propositions 11 and 12
$x\in R(A)$, since $R_c(A)=R(A)$ by the conditions of this lemma.
The algebra $A$ is semi-simple, consequently, $J=(0)$. \par For a
right ideal or a two-sided ideal the proof is analogous.

\par {\bf 18. Lemma.} {\it If $A$ is a Banach annihilator semi-simple algebra
over a field $F\in E_n$ with $R_c(A)=R(A)$ and $J_r$ is a right
minimal ideal in $A$, then a closed two-sided ideal $Y=Y(J_r)$
generated by $J_r$ is minimal and closed in $A$.}
\par {\bf Proof.} If $X$ is a closed two-sided ideal contained
in $Y$, then $J_r\cap X$ is a right ideal contained in $J_r$,
consequently, either $J_r\cap X=J_r$ or $J_r\cap X=(0)$, since $J_r$
is minimal. If $J_r\cap X=J_r$, then $Y\subset X$, hence $Y=X$.
\par If $J_r\cap X=(0)$, then $J_rX\subset J_r\cap X=(0)$,
consequently, $J_r\subset {\sf L}(A,X)$. Then ${\sf L}(A,X)$ is the
closed two-sided ideal, consequently, $Y\subset {\sf L}(A,X)$.
Therefore, $X\subset {\sf L}(A,X)$ and consequently, $X^2=(0)$.
Applying Lemma 17 we get that $X=(0)$. \par Thus $Y$ is minimal.

\par {\bf 19. Theorem.} {\it Let $A$ be a Banach annihilator semi-simple algebra
over a field $F\in E_n$ with $R_c(A)=R(A)$. Then the sum of all left
(or right) ideals of $A$ is dense in $A$.}
\par {\bf Proof.} Suppose that $U$ is a sum of all minimal right
ideals and $\bar U$ is its closure in $A$. If $\bar{U}\ne A$, then
$\bar U$ is the closed right ideal in $A$, consequently, a nonzero
element $y$ in $A$ exists such that $y\bar{U}=(0)$. This implies
that $y$ belongs to all left annihilators of all minimal right
ideals and hence it belongs to the intersection $V$ of all maximal
left regular ideals. In view of Proposition 3.8 one gets that this
intersection is $R_c(A)$. By the conditions of this theorem
$R_c(A)=R(A)$, hence $V$ is zero, since $A$ is semi-simple. Thus
$y=0$ providing the contradiction. Thus $\bar{U}=A$.

\par {\bf 20. Proposition.} {\it Let conditions of Theorem 19 be
fulfilled and let $J$ be a right ideal in $A$. Then $J$ contains a
minimal right ideal and an irreducible idempotent $s$.}
\par {\bf Proof.} Suppose that $J$ does not contain a minimal right
ideal and $sA$ is a minimal right ideal for some irreducible
idempotent $s$ in $A$. This implies that $J\cap (sA)=(0)$. Hence for
each $a\in A$ either $asA=(0)$ or $asA$ is also a minimal right
ideal, consequently, $(asA)\cap J=(0)$ for all $a\in A$ and hence
$(as)\cap J=(ass)\cap J\subset (asA)\cap J=(0)$ for all $a\in A$.
Thus $(aS)\cap J=(0)$. Therefore $JAs=(0)$, since $JAs\subset
(As)\cap J$. This means that $JAs=(0)$ for all minimal left ideals
$As$. In view of Theorem 19 $JA=(0)$, consequently, $J=(0)$.

\par {\bf 21. Proposition.} {\it If conditions of Theorem 19 are satisfied and
$s$ is an irreducible idempotent in $A$, then $sA$ and $As$ are
minimal right and left ideals correspondingly.}
\par {\bf Proof.} Suppose that $sA$ is not minimal. By virtue of
Proposition 20 it contains a minimal right ideal $rA$ such that
$rA\ne sA$, $rA\subset sA$. Then an element $a\in A$ exists such
that $r=sa$, consequently, $rs=sas\in rA$. This implies that $t$ is
a nonzero idempotent contained in $rA$ such that the element $t=rs$
satisfies the equalities $st=ts=t$ and $s-t$ is also a nonzero
idempotent providing the contradiction, since $s=t+(s-t)$ and
$t(s-t)=(s-t)t=0$, but $s$ is irreducible by the conditions of this
proposition. Thus $sA$ is minimal.

\par {\bf 22. Proposition.} {\it If conditions of Theorem 19 are satisfied
and $J$ is a closed two-sided ideal in $A$, then ${\sf L}(A,J)={\sf
R}(A,J)$ and $J+{\sf R}(A,J)$ is dense in $A$.}
\par {\bf Proof.} In view of Lemma 17 $J\cap {\sf R}(A,J)=(0)$, since
$J\cap {\sf R}(A,J)=:V$ is the right ideal possessing the property
$V^2=V$. Therefore, ${\sf R}(A,J)J=(0)$ and hence ${\sf
R}(A,J)\subset {\sf L}(A,J)$. Similarly ${\sf L}(A,J)\subset {\sf
R}(A,J)$, consequently, ${\sf L}(A,J)={\sf R}(A,J)$.
\par If $J+{\sf R}(A,J)$ would be not dense in $A$, then its closure should be a
proper ideal in $A$, consequently, a nonzero element $x$ in $A$
exists such that $(J+{\sf R}(A,J))x=(0)$. Therefore $J(\alpha
x+xy)=(0)$ and ${\sf R}(A,J)(\alpha x+xy)=(0)$ for each $y\in A$ and
$\alpha \in F$, hence $(\alpha x+xy)\in {\sf R}(A,J)$ and
consequently, $(\alpha x+xy)^2=0$ for each $y\in A$ and $\alpha \in
F$. But in the semi-simple algebra $A$ with $R_c(A)=R(A)$ this is
impossible for $x\ne 0$.

\par {\bf 23. Proposition.} {\it If conditions of Theorem 19 are met and $J$ is a minimal
closed two-sided ideal in $A$, then $J$ is an annihilator algebra
with $R_c(J)=R(J)$. If in addition $A$ is dual, then $J$ is also
dual.}
\par {\bf Proof.} If $x\in J$ and $Jx=(0)$, then $x=0$, since
$J\cap {\sf R}(A,J)=(0)$ due to Proposition 22. Analogously if
$xJ=(0)$ and $x\in J$, then $x=0$. Thus ${\sf L}(A,J)={\sf
R}(A,J)=(0)$.
\par If $V_l$ is a closed left ideal in $J$, then
$(J+{\sf L}(A,J))V_l=JV_l\subset V_l$, hence $AV_l\subset V_l$,
since $J+{\sf L}(A,J)$ is dense in $A$ by Proposition 22. Thus $V_l$
is the closed left ideal in $A$.
\par Put $H_l=V_l+{\sf R}(A,J)$. Then either $H_l$ is dense in $A$ or
${\sf R}(A, H_l)\ne (0)$. From Lemmas 17, 18 and Proposition 20 one
gets $J\cap {\sf R}(A,H_l)\ne (0)$ and hence $J\cap {\sf
R}(A,V_l)\ne (0)$. Analogously $J\cap {\sf L}(A,V_r)\ne (0)$ for a
closed right ideal $V_r$ in $J$.
\par Suppose now that the algebra $A$ is dual. In view of Lemma 17
and Proposition 22 if $x\in J$ and $[{\sf L}(A,V_r)\cap J]x=(0)$,
then $x\in {\sf R}(A,{\sf L}(A,V_r)\cap J)=cl_A({\sf R}({\sf
L}(A,V_r))+{\sf R}(A,J))=cl_A(V_r+{\sf R}(A,J))=cl_A(V_r+{\sf
L}(A,J))$. Then $(V_r+{\sf L}(A,J))J=V_rJ\subset V_r$, since $V_r$
is a right ideal in $J$, consequently, $cl_A(V_r+{\sf
L}(A,J))J\subset V_r$, hence $xJ\subset V_r$ and consequently, ${\sf
L}(A,V_r)xJ\subset {\sf L}(A,V_r)V_r=(0)$. On the other hand, ${\sf
L}(A,V_r)x{\sf R}(A,J)=(0)$, since $x\in J$, consequently, ${\sf
L}(A,V_r)x(J+{\sf R}(A,J))=(0)$. We have that $J+{\sf R}(A,J)$ is
dense in $A$ due to Proposition 22, hence ${\sf L}(A,V_r)xA=(0)$ and
consequently, ${\sf L}(A,V_r)x\subset {\sf L}(A,A)=(0)$. From the
duality of $A$ it follows that $x\in V_r$. Therefore, ${\sf
R}(J,{\sf L}(J,V_r))=V_r$ and similarly ${\sf L}(J,{\sf
R}(J,V_l))=V_l$. Thus $J$ is also dual.

\par {\bf 24. Theorem.} {\it Let $A$ be a Banach semi-simple annihilator
algebra over $F\in E_n$ with $R_c(A)=R(A)$. Then $A$ is the
completion of the direct sum of all its minimal closed two-sided
ideals $H_k$ each of which is a simple annihilator algebra over $F$.
Moreover, if $A$ is dual, then each $H_k$ is simple and dual.}
\par {\bf Proof.} By virtue of Proposition 20 each closed minimal
two-sided ideal $J$ in $A$ contains a minimal right ideal $V_r$,
hence $J=V_r$ according to Lemma 18. Then the closure $cl_AV_r$ is a
closed minimal two-sided ideal for each minimal right ideal $V_r$
due to the same lemma. According to Proposition 23 $cl_AV_r$ is the
annihilator algebra, which is also dual if $A$ is dual. If $H$ is a
closed two-sided ideal in $cl_AV_r$, then it is such in $A$ also.
But $cl_AV_r$ is minimal, hence the algebra $cl_AV_r$ is simple.
\par By virtue of Theorem 19 the sum of all minimal right ideals
$V_r$ is dense in $A$. Let $K$ and $M$ be two minimal closed
two-sided ideals which are different, $K\ne M$. Therefore $KM\subset
K\cap M=(0)$, since $K\cap M$ is the closed two-sided ideal
contained in minimal closed two-sided ideals $K$ and in $M$ and
different from them. If $x+y=0$ for some $x\in K$ and $y\in M$, then
$Kx=(0)$ and $My=(0)$, consequently, $(xA)^2\subset K(xA)=(0)$ and
analogously $(yA)^2=(0)$. Therefore $xA=(0)$ and $yA=(0)$, since $A$
is semi-simple, consequently, $x=0$ and $y=0$. Thus the considered
sum is direct.

\par {\bf 25. Theorem.} {\it If $A$ is a Banach simple annihilator algebra
over a field $F\in E_n$ with $R_c(A)=R(A)$, if also $p$ is an
irreducible idempotent, then $pAp=:H$ is an ultranormed division
algebra over $F$. Moreover, if $A$ and $F$ are ultranormed and $A$
is commutative, then a multiplicative ultranorm $|\cdot |_H$ on $H$
exists extending that of $F$ such that it induces a topology on $H$
not stronger than the topology inherited from $A$.}
\par {\bf Proof.} From the conditions of this proposition it follows
that $pH=Hp=H$, since $p^2=p$ and the algebra $A$ is associative.
Evidently, $H$ is the algebra over $F$, since $A$ is the algebra
over $F$. The restriction of $p$ to $H$ is the identity on $H$,
since $ps=p^2s=p(ps)$ for each $s\in A$ and hence $pr=r$ for each
$r\in H$, similarly $rp=r$ for each $r\in H$ and hence
$pr=rp=r=prp$. For each nonzero element $r$ in $H$ the set $Ar$ is a
left ideal in $A$ and $Ar \ne (0)$ due to Condition 1$(1)$. In view
of Propositions 20 and 21 $Ar\subset Ap$ and $Ap$ is a minimal left
ideal, since $p$ is the irreducible idempotent. Thus $Ar=Ap$ and
hence an element $y\in A$ exists such that $yr=p^2=p$, consequently,
$pyr=py(pr)=(pyp)r$. Therefore, $(pyp)r=(pyp)(prp)=pyprp=pyr=pp=p$,
consequently, $pyp$ is a left inverse of $r$ in $H$. Similarly $r$
has a right inverse in $H$. Thus $H$ is the division algebra such
that $F$ is isomorphic with $Fp$ and $Fp\subset H$. From the
continuity of the algebraic operations on $A$ it follows that they
are continuous on $H$. The norm on $A$ induces a norm on $H$, since
$H\subset A$. Since $H$ is the topological ring with the continuous
quasi-inverse and $H$ possesses the unit, then $H$ is with the
continuous inverse.
\par If $A$ and $F$ are ultranormed and $A$ is commutative,
then the ultranorm $|\cdot |_A$ on $A$ induces the ultranorm on $H$
and $H$ is also commutative. Therefore, $|p|_A=|p^2|_A\le |p|_A^2$
and hence $1\le |p|_A$. On the other hand, on $H$ as the field
extension of $F$ there exists a multiplicative ultranorm $|\cdot
|_H$ extending $|\cdot |_F$ that of the field $F$ (see Proposition 5
in Section VI.3.3 \cite{bourbcomalgb}, Krull's existence theorem
14.1 and Theorem 14.2 in \cite{sch1} or 3.19 in \cite{roo}). We have
that $|1|_F=1$, $1\le |p|_A$, also $p$ plays the role of the unit in
$H$, while $|bx|_A=|b|_F|x|_A$ for each $b\in F$ and $x\in A$. \par
If $A$ is not unital, we consider the algebra $A_1$ obtained from it
by adjoining the unit. The norms on $A$ and $F$ induce the norm on
$A_1=A\oplus F$. Therefore, it is sufficient to consider the case of
the unital algebra $A$. Mention that $(1-p)^2=1-p$ and $A(1-p)$ is
the ideal in $A$ such that $A=Ap+A(1-p)$ with $Ap\cap A(1-p)=(0)$.
Moreover, $Ap=pAp=Ap^2$, since $A$ is commutative. This implies that
$H$ is isomorphic with the quotient algebra $J=A/(A(1-p))$. Then the
ultranorm on $A$ induces the quotient ultranorm on $J$ such that
$|xy|_J\le |x|_J|y|_J$ and $|xyp|_J\le |xp|_J|yp|_J$ for each $x$
and $y$ in $J$, since $pxp=xp$ and $xpyp=xyp$ for each elements $x$
and $y$ in the commutative algebra $A$. At the same time,
$|xyp|_H=|pxppyp|_H=|xp|_H|yp|_H$ for each $x$ and $y$ in $A$.
\par The ultranorm $|\cdot |_A$ on $Fp$ induced from $A$ is equivalent
with the multiplicative ultranorm $|\cdot |_F$ on $F$, since $Fp$ is
isomorphic with $F$ and consequently, $|xpypz|_A=|xp|_A|yp|_A|z|_A$
for every $xp\in Fp$, $yp\in Fp$ and $z$ in $A$, since
$|xp|_{F}=|xp|_A$. Then $|xyp|_J=|xp|_J|yp|_J$ if $xp\in Fp$ and
$yp\in Fp$, where $x$ and $y$ are in $A$. The inequality
$|p|_J^{-1}|xp|_J\le |x|_J$ is also fulfilled for each $x\in A$.
Therefore, $H$ can be supplied with a multiplicative norm $|\cdot
|_H$ extending that of $F$ and satisfying the inequality $|x|_H\le
|x|_J$ for each $x\in H$ according to Theorems 1.15 and 1.16
\cite{escanelb}.

\par {\bf 26. Proposition.} {\it Suppose that $A$ is a
Banach simple annihilator algebra over a field $F$, $F\in E_n$, also
$R_c(A)=R(A)$. Then a maximal family of orthogonal irreducible
idempotents $ \{ w_j: j \in J \} $ exists such that $\sum _jAw_j$
and $\sum _jw_jA$ are dense in $A$.}
\par {\bf Proof.} In view of Proposition 20 there are irreducible
idempotents $w_j$ in $A$. Each right ideal $B$ in $A$ contains a
minimal right ideal, consequently, it contains an irreducible
idempotent. By virtue of Zorn's lemma (see \cite{eng}or
\cite{kunenb}) a maximal orthogonal system $ \{ w_j: j \in J \} $ of
irreducible idempotents $w_j$ exists. Let $C=\sum _jAw_j$ be the sum
of all such left ideals. Suppose that $cl_AC\ne A$. Then $cl_AC$ is
a closed left ideal. Therefore, $A_r(cl_AC)$ is the right ideal
different from zero. This implies that $A_r(cl_AC)$ contains an
irreducible idempotent $p$ orthogonal to each $w_j$. But this is
impossible, since the family $ \{ w_j: j \in J \} $ is maximal. It
remains that $C$ is dense in $A$. Similarly $\sum _jw_jA$ is dense
in $A$.

\par {\bf 27. Proposition.} {\it Let $F$ be a field and let $\{ K_j: ~ j\in P \} $
be a family of division algebras such that $F$ is contained in the
center $Z(K_j)$ of $K_j$ for each $j\in P$, where $P$ is a set. Then
a minimal division algebra $K$ exists such that $K_j\subseteq K$ for
each $j\in P$.}
\par {\bf Proof.} Since $K_j$ is a division algebra, then its center
$Z(K_j)$ is a field. Take the tensor product $T=\bigotimes_{j\in
P}K_j$ of $K_j$ as algebras over the field $F$. Therefore, $T$ is an
algebra over $F$ so that $T$ may be noncommutative if at least one
of $K_j$ is noncommutative. For each $K_j$ a natural embedding $h_j:
K_j\hookrightarrow K$ exists. Moreover, $T$ contains the unit
element which can be identified with the unit of the field $F$.
\par For each proper left ideal $B$ in $T$ the intersection $B_j=B\cap
h_j(K_j)$ is a left ideal of $K_j$. In view of Theorem I.9.1 in
\cite{bourbalgch123b} $B_j=(0)$, since $F\subset T/B$ and the unit
is unique in the associative algebra $T/B$. Particularly, for a
maximal proper left ideal $B$ in $T$ this induces the embedding
$t\circ h_j$ of $K_j$ into the quotient algebra $T/B$ over the field
$F$ for each $j\in P$, where $t: T\to T/B$ denotes the quotient
$F$-linear mapping. \par Then equations $a_jx_j=b_j$ and
$y_ja_j=b_j$ with $a_j\ne 0$ and $b_j$ in $h_j(K_j)$ have unique
solutions $x_j$ and $y_j$ in $h_j(K_j)$ for each $j\in P$. For an
arbitrary $a\in T/B$ take an element $c\in t^{-1}(a)$. Then
$c+B=t^{-1}(a)$ and consequently, $h_j(K_j)\cap
t^{-1}(a)=h_j(K_j)\cap \{ c \} $, where $\{ c \} $ denotes the
singleton in $T$. At the same time, $qu=0$ for some $q$ and $u$ in
$T$ implies $t((q+B)(u+B))=0$ in $T/B$. Therefore, equations $ax=b$
and $ya=b$ with $a\ne 0$ and $b$ in $T/B$ have unique solutions $x$
and $y$ in $T/B$, since $B$ is the proper maximal left ideal in $T$
and $h_j(K_j)\cap h_i(K_i)=F$ for each $j\ne i$ in $P$. From Theorem
9.2 and Corollary 9.3 in \cite{neumplms51} it follows that an
embedding of $T/B$ into a unique-division algebra $L$ over $F$
exists. Taking the intersection of all such algebras $L$ one gets a
minimal unique-division algebra $K$ over $F$ containing $T/B$. Thus
the embedding of $K_j$ into the division algebra $K$ exists for each
$j\in P$.

\par {\bf 28. Theorem.} {\it If conditions of Proposition 27 are
fulfilled and each $K_j$ is a Hausdorff topological division algebra
with a topology $\tau _j$ such that \par $(i)$ $\tau _j|_{K_i\cap
K_j}=\tau _i|_{K_i\cap K_j}$ \\ for each $i$ and $j$ in $P$, then a
Hausdorff topology $\tau $ on $K$ exists such that an embedding
$h_j: K_j\to K$ is a homeomorphism of $(K_j, \tau _j)$ onto
$(h_j(K_j), \tau \cap h_j(K_j))$ for each $j\in P$. Moreover, if
each $K_j$ is ultranormed and \par $(ii)$ $|\cdot |_{K_j}|_{K_i\cap
K_j}=|\cdot |_{K_j}|_{K_i\cap K_j}$ \\ for each $i$ and $j$ in $P$,
where $|\cdot |_{K_j}$ denotes an ultranorm on $K_j$, then $K$ is
ultranormed.}
\par {\bf Proof.} Consider on the weak product $S=\prod_{j\in P}' K_j$
the box product topology, where each $s\in S$ has the form $s=(s_j:
\forall j ~ s_j\in K_j, ~ card \{ j: s_j\ne e_j \} <\aleph _0 )$,
where $e_j=1$ denotes the unit element in $K_j$. It induces the
corresponding topology $a_p$ on the tensor product $T$, where $T$ is
the quotient algebra $S/M$ of $S$ by the submodule $M$ having
elements of the form \par $(1)$ $(x)+(y)-(z)$ with $x_i+y_i=z_i$ for
one index $i\in P$ and with $x_j=y_j=z_j$ for each $j\ne i$ in $P$;
\par $(2)$ $(x)-(y)$ with $x_i=by_i$ for one index $i$ in $P$ and
$x_j=y_j$ for each $j\ne i$ in $P$ for every $b\in F$, $(x)$, $(y)$
and $(z)$ in $S$ (see also Chapter 3 in \cite{bourbalgch123b}). The
algebra $T$ is supplied with the multiplication prescribed by the
rule $(x)(y)=\bigotimes_{j\in P}x_jy_j$ for each $(x)$ and $(y)$ in
$T$. Due to condition $(i)$ for each $i$ there exists an algebraic
topological embedding of $K_i$ into $T$.
\par The algebra $K$ over $F$ is obtained as the unique-division
algebra $K$ over $F$ containing $T/B$ (see the proof of Proposition
27).
\par The algebra $T$ is unital, since $K_j$ is unital for each $j$.
There exists a neighborhood $W_j$ of $1$ in $K_j$ such that the
inversion is continuous on $W_j$ for each $j$. Take $W=\prod_{j\in
P}' W_j$, hence $W$ is a neighborhood of $1$ in $T$ such that the
inversion is continuous on $W$, since $S$ is supplied with the box
topology. Therefore, if $B$ is a left maximal ideal in $T$, then $B$
is closed in $T$, since algebraic operations on $T$ are continuous
and $T$ is with the continuous inverse on $W$.
\par Therefore, the box topology $a_p$ on $T$ induces the quotient $T_1$-topology
$b_p$ on $T/B$, since $a_p$ is the Hausdorff topology and $B$ is
closed in $T$. Consider a base $U$ of a topology $\tau $ on $K$
satisfying the conditions:
\par $(3)$ $U_x=U_0+x$ for each $x\in K$,
\par $(4)$ $U_x=xU_1=U_1x$ for each nonzero $x$ in $K$,
\par $(5)$ $U_0\cap (T/B)$ is the base of neighborhoods of zero in
the $b_p$ topology on $T/B$, where $U_x$ denotes a base of
neighborhoods of an element $x$ in $K$ such that $U_x\subset U$;
\par $(6)$ for each $E$ and $D$ in $U_0$ there exists $C\in U_0$
such that $C\subset E\cap D$;
\par $(7)$ $\bigcap_{V\in U_0}V= \{ 0  \} $;
\par $(8)$ for each $E\in U_0$ there exists $D\in U_0$ such that
$(D+D)\subset E$ and $(D+1)^2 \subset (E+1)$;
\par $(9)$ for each $E\in
U_0$ there exists $D\in U_0$ such that $-D\subset E$ and $(D+1)^{-1}
\subset (E+1)$, \par $(10)$ $U|_F$ provides the base of the $\tau
_j|_F$ topology on $F$. \par This is possible, since $F\subset K_i$
and condition $(i)$ is fulfilled for each $i$ and $j$ in $P$ and
since the $b_p$ topology on $T/B$ satisfies analogous to $(3-10)$
conditions due to Theorem 1.3.12 in \cite{arhtka4b}.
\par Each element of $K$ is
obtained from elements of $T/B$ by a finite number of algebraic
operations. Therefore, the intersection of all such bases $U$
satisfying conditions $(3)-(9)$ provides a minimal base possessing
these properties. In view of Theorem 1.3.12 in \cite{arhtka4b} this
induces a Hausdorff topology $\tau $ on $K$.
\par From the construction above it follows that $\tau \cap
(T/B)=b_p$, consequently, $\tau \cap h_j(K_j)$ is equivalent with
the topology $\tau \cap h_j(K_j)$ on $h_j(K_j)$ inherited from
$(T/B, b_p)$ for each $j\in P$, where $h_j$ is the algebraic
embedding as in subsection 27. Therefore $h_j$ is the homeomorphism
of $(K_j, \tau _j)$ onto $(h_j(K_j), \tau \cap h_j(K_j))$ for each
$j\in P$.
\par In particular, if $K_j$ is ultranormed for each $j$, then
$T$ is ultranormed by $|x|=\sup_j |x_j|_{K_j}$, where $|hx|=|h||x|$
for each $h\in F$ and $x\in T$. Such ultranormed topology is not
stronger than the $a_p$ topology. By condition $(ii)$ of this
theorem ultranorms $|\cdot |_{K_j}$ and $|\cdot |_{K_i}$ on $K_i\cap
K_j$ are equivalent for each $i$ and $j$ in $P$, hence there exists
an algebraic isometric embedding of $K_i$ into $T$ for each $i$. On
the other hand, $F\subset K_i$ for each $i$. This induces the
quotient ultranorm on $T/B$ relative to which $h_j$ is continuous
for each $j$. Therefore, $U_0$ and $U_1$ on $K$ can be chosen
countable and such that $V+V\subseteq V$ for each $V\in U_0$,
$WW\subseteq W$ for each $W\in U_1$. Thus $K$ is ultranormable:
$|x+y|_K\le \max (|x|_K, |y|_K)$ and $|xy|_K\le |x|_K |y|_K$ for
each $x$ and $y$ in $K$.

\par {\bf 29. Corollary.} {\it If conditions of Theorem 28 are
satisfied, then a completion $\tilde{K}$ of $K$ relative to a left
uniformity $l_{\tau }$ induced by $\tau $ exists such that
$\tilde{K}$ is a division algebra. Moreover, if Condition 28$(ii)$
is fulfilled, then $\tilde{K}$ is the Banach division algebra.}
\par {\bf Proof.} Consider on the multiplicative group $K^*$ of nonzero
elements of $K$ the left uniformity $l_{\tau }$ induced by $\tau $.
In view of \S 8.1.17 and Theorem 8.3.10 in \cite{eng} and conditions
28$(3,4)$ the completion $\tilde{K}$ of $K$ relative to $l_{\tau }$
is the unique-division algebra over $F$. \par If in addition
Condition 28$(ii)$ is fulfilled, then we take $\tilde{K}$ as the
completion of $K$ relative to its ultranorm.

\par {\bf 30. Theorem.} {\it Let $A$ be a simple annihilator Banach
algebra over an ultranormed field $F$, $F\in E_n$, with
$R_c(A)=R(A)$, and let $ \{ w_j \} $ be a maximal orthogonal system
of irreducible idempotents in it. Then an ultranormed Banach
division algebra $G$ exists such that $w_jAw_j\subset G$ for each
irreducible idempotent $w_j$ in $A$, also $\sum_{i,j}w_iAw_j$ is
dense in $A$.}

\par {\bf Proof.} Suppose that $ \{ w_j \} $ is a maximal orthogonal
system of irreducible idempotents in the algebra $A$. For a chosen
idempotent $w_i$ one gets the two-sided non nil ideal $Aw_iA$. Since
$A$ is simple, then $cl_A(Aw_iA)=A$, where $cl_AS$ denotes the
closure of a subset $S$ in $A$. Therefore, $w_jAw_iAw_j\ne (0)$ for
each $j$. Moreover, $w_jAw_iAw_j\subseteq w_jAw_j=G_jw_j$, where
$G_j$ is a division algebra over $F$ according to Theorem 25.

\par Consider the algebra $A_j:=A_{G_j}=A\hat{\otimes }_FG_j$ obtained from $A$ by
extension. For each $j$ elements $x_j$ and $y_j$ in $A_j$ exist such
that $w_jx_jw_iy_jw_j=w_j$. Put $w_{j,i}=w_jx_jw_i$ and
$w_{i,j}=w_iy_jw_j$ and $w_{j,k}=w_{j,i}w_{i,k}$. Therefore,
$w_{j,i}$ and $w_{i,j}$ belong to $A_j$ and $w_{j,j}=w_j$, since
$w_i^2=w_i$. Then one infers that
\par $w_{j_1,k_1}w_{k_1,k_2}=w_{j_1}x_{j_1}w_iw_iy_{k_1}w_{k_1}w_{k_1}x_{k_1}w_iw_iy_{k_2}w_{k_2}$
\par $=w_{j_1}x_{j_1}w_i(y_{k_1}w_{k_1}x_{k_1}w_i)y_{k_2}w_{k_2}$
and
\par
$w_{j_1,k_2}=w_{j_1,i}w_{i,k_2}=w_{j_1}x_{j_1}w_iw_iy_{k_2}w_{k_2}$\par
$= w_{j_1}x_{j_1}w_iy_{k_2}w_{k_2}$.
\par Mention that this construction implies $w_{i,j}w_{j,i}\in
w_iAw_i=G_iw_i$ and consequently, $w_{i,j}w_{j,i}=bw_i$ for a scalar
$b=b_{i,j}\in G_i$. The multiplication of both sides of the latter
equality on the left by $w_{j,i}$ and on the right by $w_{i,j}$
leads to
\par $w_{j,i}w_{i,j}w_{j,i}w_{i,j}=w_j^2=bw_{j,i}w_{i,j}=bw_j^2$, \\
consequently, $w_j=bw_j$ and hence $b=1$. Thus
\par
b$w_{j_1,k_1}w_{k_1,k_2}=w_{j_1,i}w_{i,k_1}w_{k_1,i}w_{i,k_2}=w_{j_1,i}w_iw_{i,k_2}=w_{j_1,k_2}$,
\\ since $w_{j,i}w_i=w_{j,i}$, $w_iw_{i,j}=w_{i,j}$. Then
\par $w_{j_1,k_1}w_{j_2,k_2}=w_{j_1,i}w_{i,k_1}w_{j_2,i}w_{i,k_2}$
and
\par $w_{i,k_1}w_{j_2,i}=w_i y_{k_1}w_{k_1}w_{j_2}x_{k_2}w_i$,
consequently, \par $(1)$ $w_{j_1,k_1}w_{j_2,k_2}=0$ and
\par $(2)$ $w_{k_1}w_{j_2}=0$ for each $k_1\ne j_2$, also
\par $(3)$ $w_{j_1,k_1}w_{k_1,k_2}=w_{j_1,k_2}$ for every $j_1, k_1, k_2$.
\par Thus the set $w_jAw_k$ is composed of elements which are
multiples of the element $w_{j,k}$, consequently, $w_jxw_{k,j}\in
w_jAw_j=G_jw_j$, where the division algebra $G_j$ is over the field
$F$ according to Theorem 25. Therefore, a scalar $b\in G_j$ exists
such that $w_jxw_{k,j}=bw_j$. Multiplying on the right by $w_{j,k}$
and using $(3)$ we infer that $w_jxw_k=bw_{j,k}$, where
$b=b(j,k,x)\in G_j$. This implies that $\sum_{j,k}w_jAw_k=:B\subset
A$, where $B$ is an algebra over $F$.
\par By virtue of Theorem 28 and Corollary 29 an ultranormed Banach division algebra $G$
exists such that $G_j\subset G$ for each $j$, since $G_j$ is the
algebra over $F$ for each $j$, also since $A$ is the ultranormed
Banach algebra.
\par We put $A_G=A{\hat \otimes }_FG$, that is $A_G$ is the right
$G$-module and the algebra over $F$. Thus $\sum_jw_jAw_j=:E\subset
A_G$.
\par On the other hand, $w_jAw_j\subset A$ as the algebra over $F$
for each $j$, since $w_j\in A$ for each $j$. Mention that the sum of
all $w_jAw_k$ contains the $F$-linear span $Y$ of the set $(\sum_j
w_jA)(\sum_kAw_k)$. The multiplication and addition are continuous
on $A$, hence $Y$ is dense in the $F$-linear span $X$ of $(\sum_j
w_jA)A$, since $\sum_kAw_k$ is dense in $A$. In its turn $X$ is
dense in the $F$-linear span $V$ of $A^2$, since $\sum_j w_jA$ is
dense in $A$. Therefore, $E$ is dense in $A$, since $V$ is the
two-sided ideal in $A$ which is necessarily dense in $A$.

\par {\bf 31. Definition.} Let $X$ be a Banach space
over an ultranormed field $F$, $F\in E_n$, such that $X$ also has
the structure of a right $G$-module, where $G$ is a division algebra
over $F$. An operator $s\in L(X,X)$ will be called (right) quasi
finite dimensional if its range $s(X)$ is contained in a finite
direct sum $x_1G\oplus ...\oplus x_nG$ embedded into $X$ and such
that $s$ is right $G$-linear, that is $s(xb)=(sx)b$ for each $x\in
X$ and $b\in G$, where $x_1$,...,$x_n$ are nonzero vectors belonging
to $X$.

\par {\bf 32. Theorem.} {\it Let $A$ be a simple annihilator Banach algebra
with $R_c(A)=R(A)$ over an ultranormed field $F$, $F\in E_n$. Then
an ultranormed Banach division algebra $G$ exists such that
$A_G:=A\hat{\otimes }_FG$ has an embedding $T$ into the algebra
$L(X,X)$, where $X$ is a Banach space over $F$ and a right
$G$-module, such that
\par $(1)$ $T(A_G)$ contains all (right) quasi finite dimensional operators
so that $T(A_G)$ is a Banach subalgebra in $L(X,X)$ and
\par $(2)$ a dense subalgebra $B$ in $A_G$ exists whose image
$T(B)$ consists of quasi finite dimensional operators.}
\par {\bf Proof.} Let an ultranormed Banach division algebra $G$ be provided by Theorem 30.
Then $w_iAw_i=G_i\subseteq G$ for each $i$ and hence $w_iA_Gw_i =G$,
since $G$ is the division algebra over $F$.
\par Denote for short $A_G$ by $A$. Then $cl_A (Aw_iAx)\ne (0)$ for each
$x\ne 0$, consequently, $Aw_iAx\ne (0)$ and hence \par $(1)$
$w_iAx\ne (0)$ for each $x\ne 0$ in $A$.
\par Next we consider a left regular representation of the algebra
$A$ by operators ${\bf L}_x$ for each $x\in A$, where ${\bf
L}_xy:=xy$ for each $y\in A$. From property $(1)$ it follows that
the left regular representation $A\ni x\mapsto {\bf L}_x$ is the
$F$-linear isomorphism. On the other hand, $G_j\subset A$ and
$G_j\subset G$ and $G_jG=GG_j=G$ for each $j$. The operator ${\bf
L}_x$ is right $G$-linear for each $x\in A$, that is ${\bf
L}_x(yb)=({\bf L}_xy)b$ for each $y\in A$ and $b\in G$, since $A$
and $G$ are associative algebras over the field $F$, also $A$ has
the structure of the right $G$-module.
\par In view of Formulas 30$(2,3)$ the operator ${\bf L}_{w_{k,i}}$ maps
the one dimensional over $G$ right module $w_iAw_j$ into $w_kAw_j$,
also ${\bf L}_{w_{j,k}}w_lAw_i =(0)$ for $k\ne l$. Since the sum
$\sum_j w_iAw_j$ is dense in $w_iA$, then ${\bf L}_{w_{k,i}}w_iAw_j$
is the one-dimensional over $G$ right module. Therefore, the
operator ${\bf L}_x$ is quasi finite dimensional for each $x$ in
$B:=\sum_{j,k}w_jAw_k$.
\par Suppose now that $V$ is a one dimensional operator in $w_iA$
over $G$ and $b\in w_iA$ is an element such that $bA\ne 0$.
Therefore, $w_iA=Gb\oplus N(V)$, where $N(V) := \{ x\in w_iA: Vx=0
\} $.
\par Suppose that ${\sf L}(A,N(V))=(0)$. This implies that
${\sf L}(A,N(V)A)={\sf L}(A,N(V))=(0)$, since the closed right
$G$-module $M_{N(V)A}$ generated by $N(V)A$ has the natural
embedding $\psi $ into $A$ and $\psi M_{N(V)A}$ is a right ideal in
$A$. Therefore $A=\psi M_{N(V)A}$. On the other hand,
$N(V)=w_iN(V)$, consequently, $M_{w_iN(V)w_iA}=w_iA$. Then the
identity $w_iAw_i=w_iG$ would imply that $N(V)=w_iA$ providing the
contradiction. This implies that ${\sf L}(A,N(V))\ne (0)$.
\par Take now $x\ne 0$ in ${\sf L}(A,N(V))$. Let
$xb=0$, hence $w_iAx=M_{Gb\oplus N(V)}x=(0)$ contradicting Property
$(1)$, consequently, $xb\ne 0$ and hence $xbA$ is a non null right
ideal in $A$. Then $xbA=w_iA$, since $w_iA$ is the minimal right
ideal in $A$ and $xbA\subseteq w_iA$. Thus an element $y\in A$
exists fulfilling the condition $yxb=Vb$. This implies that the
operators $V$ and ${\bf L}_{yx}$ coincide on $Gb$. Mention that
${\bf L}_{yx}N(V)=(0)=VN(V)$, since $x\in {\sf L}(A,N(V))$, hence
${\bf L}_{yx}=V$. Thus all right $G$-linear one dimensional over $G$
operators are among ${\bf L}_x$, where $x\in A$.
\par Assume that $A$ is a Banach algebra, then $G$ provided by
Theorem 30 is also a Banach division algebra over $F$. By the
continuity of the multiplication in $A$, one gets that $w_iA$ is a
closed $F$-linear subspace in $A$, consequently, $|{\bf L}_x| \le
|x|$ for each $x\in A$, since $|{\bf L}_xy|=|xy|\le |x| |y|$ for
each $y\in A$. Therefore $A\ni x\mapsto {\bf L}_x$ is the continuous
isomorphism into $L(X,X)$ and each ${\bf L}_x$ is the limit relative
to the operator norm topology of quasi finite dimensional operators
${\bf L}_{x_n}$ with $x_n \in B$ for each $n\in \bf N$.

\section{$B^*$-algebras.}
\par {\bf 1. Definition.} Let $A$ be an ultranormed algebra
over $F\in E_n$ satisfying the following conditions:
\par $(1)$ $A$ is a Banach $*$-algebra and
\par $(2)$ there exists a bilinear form
$(\cdot , \cdot ): X^2\to F$ such that $|(x,y)|\le q |x| |y|$ for
all $x$ and $y$ in $A$, where $0<q<\infty $ is a constant
independent of $x$ and $y$, \par $(3)$ $(x,y)=(y,x)$ and
$(x,y)=(x^*,y^*)$ for each $x$ and $y$ in $A$,
\par $(4)$ if $(x,y)=0$ for each $y\in A$, then $x=0$;
\par $(5)$ $(xy,z)=(x,zy^*)$ for every $x$, $y$ and $z$ in $A$,
\par $(6)$ $xx^*\ne 0$ for each nonzero element $x\in A\setminus (0)$.
\par Then we call $A$ a $B^*$-algebra.

\par {\bf 2. Lemma.} {\it For a $*$-subalgebra $A$ of $L(X,X)$
with $X=c_0({\bf N},F)$, $F\in E_n$, a bilinear form $(\cdot , \cdot
)$ satisfying conditions 1$(2,3,5)$ exists.}

\par {\bf Proof.} We put $(x,y)=Tr (x^*Sy)$, where $S$ is a
marked compact operator such that $S^*=S$, $S\in L_c(X,X)$,
$X=c_0({\bf N},F)$, the trace $Tr (C) =\sum_jC_{j,j}$ is defined for
each compact operator $C\in L_c(X,X)$. In view of Theorem 2.4 $Tr
(x^*Sy)$ exists for each $x$ and $y$ in $L(X,X)$. Since $|Tr (C)|
\le |C|$ and $|x^*Sy|\le |x||S||y|$, then condition 1$(2)$ is valid.
From $Tr (C^*) = (Tr (C))^*=Tr (C)$ for each  $C\in L_c(X,X)$ and
$(x^*Sy)^*=y^*Sx$ property 1$(3)$ follows, since $t^*=t$ for each
$t\in F$. Then using the identity $Tr (CD)=Tr
(DC)=\sum_{k,j}C_{k,j}D_{j,k}$ for each $C\in L_c(X,X)$ and $D\in
L(X,X)$ we deduce that $(xy,z)=Tr (y^*x^*Sz) = Tr
(x^*Szy^*)=(x,zy^*)$ for every $x$, $y$ and $z$ in $A$, since
$(xy)^*=y^*x^*$.

\par {\bf 3. Lemma.} {\it If conditions of Lemma 2 are satisfied and
$L_c(X,X)\subseteq A$, then conditions 1$(4,6)$ are also valid.}
\par {\bf Proof.} Choose $S\in L_c(X,X)$ for which the decomposition
$S=T^{-1}YT$ is such that $T: X\to X$ is an automorphism of the
Banach space $X$ and $S^*=S$, also $Ye_j=Y_{j,j}e_j$ with
$Y_{j,j}\ne 0$ for each $j$, while $Y_{i,j}=0$ for each $i\ne j$,
where $ \{ e_k: k\in {\bf N} \} $ is the standard basis of $X$. Then
we get property 1$(4)$, since $Tr (x^*Sy)\in F$. \par On the other
hand, $(ax)(ax)^*=a(xx^*)a^*$ and $(xx^*)^*=xx^*$. Therefore,
considering $a\in A$ of the form $a=\sum_{k,j} a_{j,k}E_{j,k}$ with
$a_{j,k}\in F$ one finds coefficients $a_{j,k}$ such that
$(ax)(ax)^*\ne 0$, since $E_{j,k}\in L_c(X,X)$ for each $j, k$ and
$L_c(X,X)\subseteq A$, where $E_{j,k}={e'}_j\otimes e_k$,
${e'}_j=\theta (e_j)$ for each $j$ (see also Definition 2.5).
Mention that $(ax)(ax)^*\ne 0$ implies that $xx^*\ne 0$, since the
algebra $A$ is associative. Thus property 1$(6)$ also is fulfilled.

\par {\bf 4. Lemma.} {\it If $J_r$ and $J_l$ are proper or improper right and left ideals
in a $B^*$-algebra $A$, then ${\sf L}(A,J_r)$ and ${\sf R}(A,J_l)$
are orthogonal relative to the family of bilinear functionals $\{
(\cdot , \cdot )_a: a \in A \} $ complements of the sets $J_r^*$ and
$J_l^*$ in the Banach space $A$, where $(x,y)_a=(ax,ay)$ for every
$a$, $x$ and $y$ in $A$.}
\par {\bf Proof.} If $x\in {\sf L}(A,J_r)$, then $xJ_r=(0)$, hence
$(axJ_r,aA)=0$ for each $a\in A$ and consequently, $(ax,aAJ_r^*)=0$
by identity 1$(5)$ and inevitably $(ax,aJ_r^*)=0$. This means that
$x\in A\ominus J_r^*$ relative to $\{ (\cdot , \cdot )_a: a \in A \}
$, that is ${\sf L}(A,J_r)$ is the orthogonal complement of $J_r^*$.
Similarly ${\sf R}(A,J_l)$ is the orthogonal complement of $J_l^*$
in $A$ as the Banach space relative to the family $\{ (\cdot , \cdot
)_a: ~ a \in A \} $ of bilinear functionals.

\par {\bf 5. Proposition.} {\it Any $B^*$-algebra $A$ is dual.}
\par {\bf Proof.} If $J_r$ and $J_l$ are right and left ideals in $A$,
then by Lemma 4 ${\sf R}(A,{\sf L}(A,J_r)) = {\sf R}(A,A\ominus
J_r^*)=A\ominus (A\ominus J_r)=J_r$ and analogously ${\sf L}(A,{\sf
R}(A,J_l))=J_l$, since $A^*=A$ and $(J_r^*)^*=J_r$.

\par {\bf 6. Theorem.} {\it Any $B^*$-algebra $A$ over a
field $F\in E_n$ with $R_c(A)=R(A)$ is representable as the direct
sum of its two-sided minimal closed ideals which are simple
$B^*$-algebras and pairwise orthogonal relative to the family of
bilinear functionals $\{ (\cdot , \cdot )_a: ~ a \in A \} $.}
\par {\bf Proof.} By virtute of Theorem 3.24 and Proposition 5 the algebra
$A$ is the completion (relative to the ultranorm) of the direct sum
of its minimal closed two-sided ideals which are simple dual
subalgebras. Consider a two-sided minimal closed non null ideal $J$
in $A$. The involution mapping $x\mapsto {\cal I}x=x^*$ provides
from it the minimal closed two-sided ideal $J^*$ due to Condition
2.5$(1)$.
\par Suppose that $J^*\ne J$, then $JJ^*=(0)$, since the ideal $J$
is minimal. From $aJ\subset J$ and $Ja\subset J$ for each $a\in A$
we deduce that $AJJ^*A=(0)$. Together with condition 1$(6)$ imposed
on the $B^*$-algebra this would imply that $x=0$ for each $x\in J$
contradicting $J\ne (0)$. Thus $J^*=J$.
\par Mention that properties 1$(1-3)$ and 1$(5)$ for $J$ are inherited from that of $A$.
Then condition 1$(6)$ on $A$ implies that $J^2\ne (0)$, since
$J^*=J$ and $AJ\subseteq J$, also $JA\subseteq J$. But $J$ is
minimal, hence $J^2=J$. Therefore, property 1$(4)$ on $J$ follows
from that of on $A$ and 1$(5)$ and $J^2=J$, since for each $u\in J$
there exists $x$ and $y$ in $J$ with $u=xy$ and
$(u,z)=(xy,z)=(x,zy^*)$ for all $z\in A$, also since $zy^*\in J$.
Then for each $y\in J\setminus (0)$ an element $x\in J\setminus (0)$
exists such that $xy\ne 0$, hence $u=xy\in J\setminus (0)$. Then we
have that $uu^*\ne 0$ by 1$(6)$ on $A$. Hence $(xy)(xy)^*\ne 0$,
consequently, $yy^*\ne 0$, since the algebra $A$ is associative and
$x(yy^*)x^*\ne 0$. Therefore property 1$(6)$ on $J$ is valid. Thus
$J$ is the $B^*$-algebra.
\par If $J$ and $S$ are two distinct minimal closed two-sided ideals
in $A$, then $JS=(0)$. From Lemma 4 it follows that $S\subset {\sf
R}(A,J)=A\ominus J^*=A\ominus J$. Thus these ideals $J$ and $S$ are
orthogonal relative to the family $ \{ (\cdot , \cdot )_a: a\in A \}
$ of bilinear functionals.
\par Using condition 1$(4)$ and Lemma 4 we infer that $A$ is the
direct sum of its two-sided minimal closed ideals.

\par {\bf 7. Theorem.} {\it Let $A$ be a simple $B^*$-algebra over a
field $F\in E_n$ with $R_c(A)=R(A)$ and let a division algebra $G$
be provided by Theorem 3.30. Then the following conditions are
equivalent:
\par $(1)$ $A_G$ is finite dimensional over $G$; \par $(2)$ $A_G$ is
unital;
\par $(3)$ the center $Z(A_G)$ of $A_G$ is non null.}
\par {\bf Proof.} Let $ \{ w_j: ~ j\in \Lambda \} $ be a maximal
system of irreducible idempotents provided by Theorem 3.30.
\par $(1)\Rightarrow (2)$. If $A_G$ is finite dimensional over $G$,
then according to Theorem 6 a maximal system $\{ w_j: ~ j\in \Lambda
\} $ of irreducible idempotents is finite, that is $card (\Lambda
)<\aleph _0$. Then their sum $w = \sum_{j\in \Lambda } w_j$ is the
idempotent fulfilling the condition $x=\sum_{j\in \Lambda } xw_j=xw$
and $x=\sum_{j\in \Lambda } w_jx=wx$. Thus $w$ is the unit in $A_G$.
\par $(2)\Rightarrow (3)$. If $A_G$ contains a unit $w$, then $Z(A_G)$
contains $w$, consequently, $Z(A_G)$ is non null.
\par $(3)\Rightarrow (1)$. Let $Z(A_G)\ne (0)$ and $x$ be a non zero element of
$Z(A_G)$, $x\ne 0$. In view of Theorem 3.30
$xw_j=(xw_j)w_j=w_jxw_j=w_j^2xw_j$, hence $xw_j=b_jw_j=w_jb_jw_j$,
where $b_j\in G$. Thus $(b_jw_j)w_j=w_j(b_jw_j)$. Therefore $x=
\sum_j xw_j = \sum_j b_j w_j$ and hence $b_jw_{j,k} = b_j w_j
w_{j,k} = x w_j w_{j,k} = x w_{j,k} = w_{j,k} x = w_{j,k}w_kx=
w_{j,k} xw_k = w_{j,k} b_kw_k$. Similarly
$b_kw_{k,j}=w_{k,j}b_jw_j$, consequently, $b_jw_{j,k}w_{k,j} =
b_jw_j =w_{j,k} b_kw_kw_{k,j} = w_{j,k} b_kw_{k,j}$ and hence
$\sum_j b_j w_j= b_kw_k + \sum_{j, j\ne k} w_{j,k} b_kw_{k,j}
=\sum_j w_{j,k}b_kw_{k,j} $.
\par Mention that $w_jA_Gw_j=Gw_j$ for each $j$, where $w_j$ plays the role of the unit in $Gw_j$.
Then \par $Gw_j\supseteq
w_j(w_{j,k}A_Gw_{k,j})w_j=w_{j,k}A_Gw_{k,j}$ \par $=
w_{j,k}(w_kA_Gw_k)w_{k,j}\supseteq
w_{j,k}(w_{k,j}A_Gw_{j,k})w_{k,j}= w_jA_Gw_j=Gw_j$ \\ for each $j$
and $k$, hence $Gw_k\ni b \mapsto w_{j,k}bw_{k,j}\in Gw_j$ is the
isomorphism of ultranormed algebras $Gw_j$ with $Gw_k$ for each $j$
and $k$.
\par  Therefore the sum $\sum_j w_{j,k}b_kw_{k,j} =
\sum_j w_{j,k}w_kb_kw_kw_{k,j}$ may converge only if it is finite.
Thus the algebra $A_G$ is finite dimensional over $G$.

\par {\bf 8. Notation.} For the Banach space $X=c_0(\alpha ,F)$ over a field
$F\in E_n$ by $L_d(X,X)$ is denoted the space of all bounded
$F$-linear operators $U: X\to X$ satisfying the condition
$\lim_{j,k} U_{j,k}=0$, that is for each $t>0$ a finite subset
$\gamma $ in a set $\alpha $ exists such that $|U_{j,k}|<t$ for each
$j$ and $k$ with either $j\in \alpha \setminus \gamma $ or $k\in
\alpha \setminus \gamma $. Then for a division algebra $H$ over $F$
and a Banach two-sided $H$-module $X_H=c_0(\alpha ,H)$ by
$L_{r,d}(X_H,X_H)$ we denote the Banach right $H$-module of all
bounded $F$-linear right $H$-linear operators $C$ from $X_H$ into
$X_H$ of the class $L_d$, that is $C(xb)=(Cx)b$ for each $x\in X_H$
and $b\in H$.

\par {\bf 9. Theorem.} {\it Let $A$ be a spherically complete
simple $B^*$-algebra over a spherically complete field $F\in E_n$
with $R_c(A)=R(A)$. Let also $G$ be a division algebra provided by
Theorem 3.30 such that $s^{1/2}\in G$ for each $s\in G$, also
$G\subset A$ and $G^*=G$. Then a Banach two-sided $G$-module $X_G$
exist such that $A$ and $L_{r,d}(X_G,X_G)$ are isomorphic as the
Banach right $G$-modules and as $F$-algebras.}

\par {\bf Proof.} By the conditions of this theorem a division algebra $G$ is such that
$wAw\subset Gw$ for each irreducible idempotent $w$ in $A$. Put
$H=G\cap G^*$. From $G=G^*$ it follows that $H=G$. If $b\in H$, then
$b^{1/2}\in G$ and $(b^{1/2})^*=(b^*)^{1/2}\in G$, since $H^*=H$,
consequently, $b^{1/2}\in H$. \par For each irreducible idempotent
$w$ described in the proof of Theorem 3.32 $ww^*\ne 0$, since $A$ is
the $B^*$-algebra over $F$. Then $(ww^*)(ww^*)^*\ne 0$, hence
$ww^*ww^*\ne 0$ and consequently, $w^*ww^*\ne 0$ implying that
$ww^*w\ne 0$, since $(w^*ww^*)^*= ww^*w$ and $c^{**}=c$ for each
$c\in A$. Therefore $w^*w\ne 0$ also.
\par Since $w$ is the irreducible idempotent and $A^*=A$, then $w^*$
is the irreducible idempotent in the $B^*$-algebra $A$. Then we
deduce that $w^*ww^*\in (w^*A_Gw^*)w^*\subseteq G^*w^*=(wG)^*$,
since $A^*=A$, consequently, an element $s\in G^*\setminus (0)$
exists such that $w^*ww^*=sw^*$, since $w^*ww^*\ne 0$. The latter
implies $w^*ww^*w=sw^*w$. But the elements $w^*ww^*w$ and $w^*w$ are
self-adjoint, hence $sw^*w=w^*ws^*$ and consequently,
\par $w^*w(s^*)^{-1}=s^{-1}w^*w$.
\par  We put $v=s^{-1}w^*w$,
hence \par $v^* = w^* w(s^*)^{-1}= s^{-1} w^* w = v$ and
\par $v^2=s^{-1}w^*ws^{-1}w^*w=
s^{-1}w^*ww^*w(s^*)^{-1}$\par $ =
(s^{-1}(sw^*w))(s^*)^{-1}=w^*w(s^*)^{-1}$ \par $=s^{-1}w^*w = v$.
\\ Thus $v$ is the self-adjoint idempotent. On the other hand,
$A_Gv=A_Gs^{-1}w^*w\subseteq A_Gw$ and $A_Gv\ne 0$ and the
idempotent $w$ is irreducible, hence the idempotent $v$ is also
irreducible, since $A_Gw$ is the non null minimal left ideal in
$A_G$.
\par Then from the proof of Theorem 3.30 it follows that
$(vA_Gv)^*=v^*A_G^*v^*=vA_Gv$ is the self-adjoint division algebra
for each such irreducible self-adjoint idempotent $v$, consequently,
$vA_Gv \subseteq Hv$. By the conditions of this theorem we have
$A=A_G$.
\par The algebra $A$ is simple, that is by the definition each its
two-sided ideal coincides with either $(0)$ or $A$.
\par Next we take a maximal orthogonal system $\{ w_j: j\in \Lambda
\} $ of self-adjoint idempotents in $A$ and for them elements
$w_{j,k}$ as in Theorem 3.30, where $\Lambda $ is a set. Hence
$w_{j,k}w_{j,k}^*\in w_jAw_j$ and $b=b_{j,k}\in H$ exists such that
$w_{j,k}w_{j,k}^*=bw_j$. Then $bw_j=w_jb^*$, since $w_j^*=w_j$ and
$(w_{j,k}w_{j,k}^*)^*=w_{j,k}w_{j,k}^*$. Moreover, $b\ne 0$, since
$w_{j,k}$ is non null and hence $w_{j,k}w_{j,k}^*$ is non null. For
$v_{j,k}=(b_{j,k})^{-1/2} ~ w_{j,k}$ we deduce that
$v_{j,k}v_{j,k}^*=w_j$, since
\par $b^{-1/2}w_{j,k}w_{j,k}^*(b^{-1/2})^*=b^{-1/2}bw_j(b^{-1/2})^*$
\par $=w_j(b^{1/2})^*(b^{-1/2})^*=w_j(b^{-1/2}b^{1/2})^*=w_j$, \\ since $A$ is
associative and $b^{-1/2}\in H$ for each non null $b$ in $H$, where
$b=b_{j,k}$. \par Thus it is possible to choose an element $w_{j,k}$
such that $w_{j,k}w_{j,k}^*=w_j$ for each $k$. Taking a marked
element $j=j_0$ and setting $w_{k,j}=w_{j,k}^*$ and
$w_{l,k}=w_{l,j}w_{j,k}$ for each $l$ and $k$ one gets
$w_{l,k}^*=w_{j,k}^*w_{l,j}^*=w_{k,j}w_{j,l}=w_{k,l}$ and
$w_{k,k}=w_k$, also $w_{k,l}w_{i,h}=\delta _{l,i} w_{k,h}$ for every
$h, i, k, l$. Thus elements $w_{l,k}$ can be chosen such that
$w_{l,k}^*=w_{k,l}$ for each $l$ and $k$.
\par If prove the statement of this theorem for the spherical completion
$\tilde H$ of $H$, then it will imply the statement of this theorem
for $H$. So the case of the spherically complete division algebra
$H$ is sufficient. Then $A$ and $H$ considered as the Banach spaces
over the spherically complete field $F$ are isomorphic with
$c_0(\alpha ,F)$ and $H$ with $c_0(\beta ,F)$ due to Theorems 5.13
and 5.16 in \cite{roo}, where $\beta \subset \alpha $. \par  From
the proof of Theorem 3.32 it follows that the sum
$B:=\sum_{j,k}w_jAw_k$ is dense in $A$. Conditions 1$(2,3,5)$ imply
that $(xy,z)=(y,x^*z)$, since $t^*=t$ for each $t\in F$. Therefore,
from properties 1$(2,3,5)$ it follows that if $j\ne h$ or $k\ne l$,
then $(w_jxw_k,w_hzw_l)=0$ for each $x$ and $z$ in $A$, since
\par $(w_jxw_k,w_hzw_l)=(w_jx,w_hzw_lw_k^*)=(w_jx,w_hz(w_lw_k))=(w_jx,0)=0$
for each $k\ne l$, also \par
$(w_hzw_l,w_jxw_k)=(zw_l,w_h^*w_jxw_k)=(zw_l,(w_hw_j)xw_k)=(zw_l,0)=0$ for each $j\ne h$. \\
Thus the set $\{ w_{j,k}: j, k \} $ is complete and
$(w_{j,k}H,w_{h,l}H)=(0)$ for each $j\ne h$ or $k\ne l$, where the
latter property is interpreted as the orthogonality. Together with
property 1$(4)$ this implies that each element $x\in A$ has the form
$x=\sum_{j,k\in \Lambda }w_{j,k} x_{j,k}$ with $\lim_{j,
k}w_{j,k}x_{j,k}=0$, since $A_H$ is the right $H$-module, also $A$
is isomorphic with $A_G$ as the $F$-algebra and the right
$G$-module, where the series may be infinite, $x_{j,k}\in H$ for
each $j, k\in \Lambda $, where $\Lambda $ denotes the corresponding
set.
\par Take the two-sided Banach $H$-module $X_H=c_0(\Lambda ,H)$ and
to each element $x\in B$ one can pose the operator $T_x$ such that
$e_j^*T_xe_k= x_{j,k} \xi _{j,k}$, where $\xi _{j,k}\in F$ and $|\xi
_{j,k}|=|w_{j,k}|$ for each $j$ and $k$ in $\Lambda $, since $|a|\in
(\Gamma _F\cup \{ 0 \} ) $ for each $a\in A$, where $B:=\sum_{j,k}
w_jAw_k$ (see above). Then $T_x\in L_{r,d}(X_H,X_H)$ and the mapping
$T: B\to L_{r,d}(X_H,X_H)$ is the isometry having the isometrical
extension $T: A\to L_{r,d}(X_H,X_H)$. The property
$w_{j,k}w_{j,k}^*=w_j\ne 0$ given above provides $|w_{j,k}|\ne 0$
for each $j$ and $k\in \Lambda $, consequently, $T$ is bijective
from $A$ onto $L_{r,d}(X_H,X_H)$, since $A$ is simple. \par For each
$S$ and $V$ in $L_{r,d}(X_H,X_H)$ one has $SV(xb)=S(Vx)b=(SVx)b$ for
each $b\in H$ and $x\in X_H$. Moreover, $|(SV)_{j,k}|\le \sup_m
|S_{j,m}| |V_{m,k}|$, consequently, $\lim_{j,k} (SV)_{j,k}=0$, that
is $SV\in L_{r,d}(X_H,X_H)$. Hence verifying other properties one
gets that $L_{r,d}(X_H,X_H)$ also has the $F$-algebra structure.
From the construction of $A_H$ it follows that $A_H$ is the
$F$-algebra, since $H$ and $A$ are $F$-algebras. Mention that
moreover, $A_H$ as the $F$-algebra is isomorphic with the Banach
$F$-algebra $L_{r,d}(X_H,X_H)$. By the conditions of this theorem
$A_H$ is isomorphic with $A$ as the $F$-algebra and the right
$H$-module.

\par {\bf 10. Theorem.} {\it Let $A$ be a spherically complete
simple $B^*$-algebra over a spherically complete field $F\in E_n$
with $R_c(A)=R(A)$ and $Z(A)=F$. Let also $G$ be a division algebra
provided by Theorem 3.30 such that $s^{1/2}\in G$ for each $s\in G$.
Then a division subalgebra $H$ of $G$ and a Banach two-sided
$H$-module $X_H$ exist such that $A_H$ and $L_{r,d}(X_H,X_H)$ are
isomorphic as the Banach right $H$-modules and as $F$-algebras.}
\par {\bf Proof.} In this case
$H=G\cap G^*$ and instead of $A$ we consider $A_H=A\hat{\otimes
}_FH$.
\par The $B^*$-algebra $A$ is simple and central, $Z(A)=F$, hence the right $H$-module $A_H$ is
simple due to Satz 5.9 in \cite{kersten} and Theorem 7 above. We
denote $A_H$ shortly by $A$ and the rest of the proof is similar to
that of Theorem 9.

\par From Theorems 6, 9 and 10 the corollary follows.

\par {\bf 11. Corollary.} {\it Suppose that $A$ is a spherically
complete $B^*$-algebra over a spherically complete field $F\in E_n$
with $R_c(A)=R(A)$ and $G$ is a division algebra given by Theorem
3.30 so that $s^{1/2}\in G$ for each $s\in G$ such that either (i)
$G\subset A$ and $G^*=G$ or (ii) $Z(A)=F$. Then a division
subalgebra $H$ in $G$ with $H^*=H$ and two-sided $H$-modules
$X_{k,H}$ exist such that $A_H$ is the direct sum of
$L_{r,d}(X_{k,H},X_{k,H})$.}

\end{document}